\documentclass[12pt]{article}
\usepackage{mathrsfs}

\usepackage{amsmath}
\usepackage{amssymb}
\usepackage{amsfonts}
\usepackage{graphicx}% 图形支持宏包
\usepackage{subfigure}
\RequirePackage[colorlinks,citecolor=blue,urlcolor=blue]{hyperref}

\newtheorem{theorem}{Theorem}
\newtheorem{definit}{Definition}
\newtheorem{lemma}{Lemma}
\newtheorem{prop}{Proposition}

\begin{document}

\title{The rate of convergence of harmonic explorer to
 $\mbox{SLE}_4$ }

\author{Shi-Yi Lan$^1$\footnote{Shi-Yi Lan was partially supported by the NSF of China (11661011) and NSF of Guangxi (2016GXNSFAA380099).
E-mail address: lanshiyi@gxun.edu.cn},
 Jin Ma$^1$\footnote{E-mail address: majingsinna@sina.com}, and Wang Zhou$^2$\footnote{Wang Zhou was partially supported by a grant
R-155-000-211-114 at the National University of Singapore. E-mail address: stazw@nus.edu.sg}\\
{}\\
{\it 1.Guangxi University for
Nationalities, 2. National University of Singapore}
}

\date{ }
\maketitle

\textbf{Abstract.} Using the estimate of the difference between the discrete harmonic function and
its corresponding continuous version we derive a rate of convergence of the Loewner driving function for
the harmonic explorer to the Brownian motion with speed 4 on the real line. 
 Based on this convergence rate, the derivative estimate
for chordal $\mbox{SLE}_4$, and the estimate of tip structure modulus for harmonic explorer paths,
we obtain an explicit power-law rate of convergence of the harmonic explorer paths to the trace of chordal $\mbox{SLE}_4$
in the supremum distance.



%\titlepage

\section{Introduction}\label{intro}

Schramm Loewner evolution with parameter $\kappa$ ($\mbox{SLE}_\kappa$),
introduced by Schramm \cite{SCH00} in 2000, is a one-parameter family of
conformally invariant random fractal curves
obtained by solving the Loewner equation with the driving term being a time-change of
one-dimensional Brownian motion. This process is intimately connected with many models in mathematical physics. It was proved that $\mbox{SLE}$s are scaling limits of a number of two-dimensional
discrete models including the Ising model \cite{SMI10,CS12,HK,CHI}, the exploration path in percolation \cite{SMI01, CN},
the loop-erased random walk (LERW) and uniform spanning trees \cite{LSW04}, the harmonic explorer \cite{SS05},
the discrete Gaussian free field \cite{SS09}. Another success of  $\mbox{SLE}$  lies in the area of probability.
For example, one can refer to \cite{LSW01A,LSW01B,LSW02} for the determination of the intersection exponents of planar Brownian motion.
As a consequence of these successes,
$\mbox{SLE}_\kappa$ processes have been intensively studied (see \cite{KN04,L05,RS05,Dub07}).
In these references, there are several different
versions of $\mbox{SLE}$s, among which the chordal $\mbox{SLE}$
and radial $\mbox{SLE}$ are the most well-known.
A chordal $\mbox{SLE}$ trace describes a
random curve evolving in a simply connected domain from one point on the boundary to another point on the boundary.
A radial $\mbox{SLE}$ trace describes a random curve evolving in a  simply connected domain from a point on the boundary to
an interior point. The behavior of the $\mbox{SLE}_\kappa$  trace depends on the real parameter $\kappa>0$.
If $\kappa\in (0,4]$, the trace is a simple curve; if $\kappa\in (4,8)$, the trace is self-touching;
and if $\kappa\in [8,\infty)$,
the trace is space-filling. For more information on chordal and radial $\mbox{SLE}_\kappa$
and related topics, one can refer to \cite{KN04,L05} and \cite{RS05}.

Despite a rapid progress on developing $\mbox{SLE}_\kappa$,
there are still several fundamental open problems. For example, one is much less aware of the speed of convergence of discrete
processes  to $\mbox{SLE}_\kappa$. One motivation for this question, besides its independent
interest, is that results of this type could lead to improved estimates of certain critical exponents \cite{sc07}.
To the best of our knowledge, the convergence rate of the planar loop-erased random walk to the radial $\mbox{SLE}_2$,
obtained by Bene\u{s}, Viklund, and Kozdron \cite{BVK}, seems to be the only result in this direction so far.
This result was improved in \cite{joh}, i.e., the rate of convergence with
respect to Hausdorff distance given by a form of non-power-law (see \cite[Throrem 8.1]{BVK}) is
improved to the one with respect to the supremum distance expressed by an explicit power-law (see \cite[Theorem 4.3]{joh}).
It is also interesting to note that the rate of convergence for Cardy's formula in site percolation has been derived by \cite{MNW,BCL}.
In this paper, we will consider the convergence rate of the harmonic explorer to  $\mbox{SLE}_4$.
The harmonic explorer is a random grid path induced by discrete harmonic
functions. In \cite{SS05} Schramm and Sheffield have proved that the harmonic explorer converges to
the chordal $\mbox{SLE}_4$ as the mesh size tends to zero. %Then it is natural to ask  whether one can estimate
%the convergence rate. In this paper, we shall give an affirmative answer to this question.

Before we state the main results of this paper, let us introduce some notation.
For any $\epsilon>0$, let $TG^\epsilon$ denote the  triangular grid of the plane $\mathbb{C}$ with the lattice
$V^\epsilon=\{\epsilon(m+ne^{i\pi/3}):m,n\in\mathbb{Z}\}$. Suppose that $D\subset\mathbb{C}$ is
a bounded Jordan domain $D\subset\mathbb{C}$
with two prescribed boundary points $u_0$ and $u_e$. Let $D^\epsilon$ be the $TG^\epsilon$
domain approximation of $D$, i.e., the largest connected component of $TG^\epsilon$
contained in $D$. Let $\hat{v}_0$ and
$\hat{v}_e$ be the midpoints of two boundary edges of $D^\epsilon$, respectively, closest to $u_0$ and $u_e$.
Then $\partial D^\epsilon$ is partitioned into two components: $L_1$, the positively oriented arc from $\hat{v}_0$ to $\hat{v}_e$
and $L_2$, the negatively oriented arc from $\hat{v}_0$ to $\hat{v}_e$. We color all vertices of $L_1\cap V^\epsilon$
white, and all vertices of $L_2\cap V^\epsilon$ black. Let $\gamma^\epsilon:[0,N]\rightarrow D^\epsilon\cup\{\hat{v}_0,\hat{v}_e\}$
be the harmonic explorer path  from $\hat{v}_0$ to $\hat{v}_e$; see Section \ref{haex} for a precise definition.
Let $\phi$ be the conformal map from $D^\epsilon$ onto the upper
half-plane $\mathbb{H}$ with $\phi(\hat{v}_0)=0$ and
$\phi(\hat{v}_e)=\infty$. Let $W^\epsilon(t)$ denote the Loewner driving function for the curve
$\tilde{\gamma}^\epsilon=\phi(\gamma^\epsilon)$ parameterized by capacity from $\infty$ in $\mathbb{H}$.
Then our first result is the following theorem on the convergence rate of
the driving function $W^\epsilon(t)$ to the Brownian motion.

\begin{theorem}\label{thm1}
Fix $T>1$. There always exists an $\epsilon_{0}>0$ depending only on $T$ with the following property.
For any $\epsilon<\epsilon_{0}$, there is a coupling of $\gamma^\epsilon$ with one dimensional standard Brownian motion $B(t),t\geq 0$,
such that
\begin{equation}\label{eq1-1}
\mathbb{P}\Big(\sup_{0\leq{t}\leq{T}}\mid{W^\epsilon(t)
-B(4t)}\mid>\epsilon^{1/12-\nu} \Big)<\epsilon^{1/12-\nu}
\end{equation}
for each fixed $0<\nu<1/12$.
\end{theorem}

Under the hypothesis that the boundary of $D$ is a $\mathcal{C}^{1+\alpha}$ curve for any $\alpha>0$,
our next result is the following theorem on the convergence rate of the harmonic
explorer path $\gamma^\epsilon$ to the chordal $\mbox{SLE}_4$ path.

\begin{theorem}\label{thm2}
Let  $\tilde{\gamma}(t)$ denote the chordal $\mbox{SLE}_4$
path in the upper-half plane $\mathbb{H}$ driven by $B(4t)$ in the coupling stated in Theorem \ref{thm1}.
For every $T>0$, there exists $\epsilon_1=\epsilon_1(T)>0$ such that if
$\epsilon<\epsilon_1$, then
\begin{equation}\label{eq1-2}\mathbb{P}\{\sup\limits_{t\in[0,T]}|\tilde{\gamma}^\epsilon(t)-\tilde{\gamma}(t)|
>\epsilon^{1/180-\mu}\}<\epsilon^{1/180-\mu}
\end{equation}
for any fixed $0<\mu<1/180$, where both curves are parameterized by the half plane capacity.  Moreover, the inequality (\ref{eq1-2})
also holds when $\tilde{\gamma}^\epsilon(t)$ and $\tilde{\gamma}(t)$ are replaced by $\gamma^\epsilon(t)$ and $\gamma(t)$
respectively, where $\gamma(t)$ is the image of $\tilde{\gamma}(t)$ in $D^\epsilon$ under the conformal map $\phi^{-1}$.
\end{theorem}

The proof strategies  are as follows. To prove Theorem \ref{thm1},
we first derive a rate of convergence for the martingale
observable of harmonic explorer, which uses an error estimate between the discrete harmonic function and
continuous harmonic function with the same boundary values, and a result of Warschawski \cite{wa}. Next,
combining the chordal Loewner equation we provide the moment estimates for increments of the
driving function $W^\epsilon(t)$. Lastly, using the Skorokhod embedding theorem
and combining estimates on the difference of two martingales and the modulus of continuity of
Brownian motion we obtain the rate of convergence of
the driving function $W^\epsilon(t)$ to the scaled Brownian motion $B(4t)$; see the inequality (\ref{eq1-1}).

In order to prove Theorem \ref{thm2}, we first establish the corresponding estimate of Loewner curves
in the deterministic setting. Consider a deterministic setting with two solutions to the chordal Loewner equation driven by
functions which are at uniform distance at most $\epsilon>0$.
If the growth of the derivative of one solution is known and the Loewner curve corresponding to  the other
solution satisfies the John-type condition, then the supremum distance between the corresponding two curves can be estimated;
see Proposition \ref{pro6-2} for a precise statement. Next, we will show
that the assumptions above are satisfied with large probability,
which results in  the speed of convergence of $\tilde{\gamma}^\epsilon$ to
$\tilde{\gamma}$ in terms of a power-law of $\epsilon$; see the inequality (\ref{eq1-2}).

Although the approaches to proving Theorem \ref{thm1} and Theorem \ref{thm2} are similar to the ones of \cite[Theorem 1.1]{BVK} and
\cite[Theorem 4.3]{joh}) respectively, there are the following essential differences.  First,
the discrete models discussed are different. The LERW is investigated in \cite{BVK} and
\cite{joh}, while we deal with the harmonic explorer process.
Hence a lot of details in our proofs are different from those in \cite{BVK} and \cite{joh}.
Secondly, the martingale observables in the case of LERW and harmonic explorer are totally different. They are discrete Green function and discrete harmonic
function respectively. In fact we apply similar techniques in \cite{sc} to obtain the convergence rate of discrete harmonic function to its continuous limit.
In addition,  the estimation of tip structure modulus
for LERW in \cite{joh} is based on the probability of crossing annulus; whereas
the estimation of tip structure modulus for harmonic explorer is via the revisiting probability estimate.
From the latter  result we get the required estimate.  Thirdly, the convergence rates obtained are different.
The convergence exponent for the driving process of LERW given in \cite{BVK} is $1/24$,
whereas the convergence exponent for the driving process of harmonic explorer is $1/12$.
Moreover, the convergence exponent for the LERW path obtained in \cite{joh} is $1/984$;
however,  the convergence exponent for harmonic explorer path is $1/180$. This shows that
the speed of convergence of the harmonic explorer to chordal $\mbox{SLE}_4$ is faster than that of
LERW to $\mbox{SLE}_2$.

This paper is organized as follows. In Section \ref{chetsm}, we introduce some notation
and review briefly some basic concepts related to the chordal $\mbox{SLE}_\kappa$, harmonic explorer process
and tip structure modulus that will be used throughout this paper. The  convergence rate
for the martingale observable is given in Section \ref{rcmo}. The moment estimates for increments of the driving
function will be presented in Section \ref{meidf}. Theorem \ref{thm1}
and Theorem \ref{thm2} will be proved in Section \ref{ecrdf} and Section \ref{ecrhep}, respectively. In Section
\ref{appen} (Appendix), we will prove a result on the tip structure modulus, that is,
the tip structure modulus of the image of a curve under a conformal map can be bounded by
the one of the curve up to a multiplicative constant under the hypothesis that the boundary of domain
is sufficiently regular.

\section{Chordal $\mbox{SLE}_\kappa$, harmonic explorer and tip structure modulus} \label{chetsm}
In this section we introduce briefly some notation that will be used throughout this paper. More information
concerning the chordal $\mbox{SLE}_\kappa$, harmonic explorer process
and tip structure modulus can be found in \cite{KN04,RS05,L05,SS05,joh}.

\subsection{Chordal $\mbox{SLE}_\kappa$}

In this subsection we briefly review the definition of chordal $\mbox{SLE}_\kappa$; see \cite{KN04,RS05,L05} for more details.

We denote by $\mathbb{H}$ the upper half plane as in Section \ref{intro}.
For fixed $T>0$, let $\gamma:[0,T]\rightarrow\overline{\mathbb{H}}$
be a continuous simple curve in $\overline{\mathbb{H}}$ which satisfies
$\gamma[0,T]\cap{\mathbb{R}}=\{\gamma(0)\}=\{0\}$, where $\mathbb{R}$ denotes the real axis.
Then for each time $t\in[0,T]$, there is a unique conformal homeomorphism
$g_{t}:H_{t}=\mathbb{H}\backslash{\gamma[0,t]}\rightarrow\mathbb{H}$
 which satisfies the so-called hydrodynamic normalization at infinity
\begin{equation}\label{eq01}
\lim_{z\rightarrow\infty}g_{t}(z)-z=0.
\end{equation}
The limit
\begin{displaymath}
\mbox{hcap}_{\infty}(\gamma[0,t])=\lim_{z\rightarrow\infty}\frac{z(g_{t}(z)-z)}{2}
\end{displaymath}
exists, which is called the half plane capacity of
$\gamma[0,t]$ from $\infty$. It is obvious that $\mbox{hcap}_{\infty}(\gamma[0,t])$ is real and monotonely increasing in $t$.
Since $\mbox{hcap}_{\infty}(\gamma[0,t])$ is also continuous in $t$, it is natural to reparametrize $\gamma$ so that
$\mbox{hcap}_{\infty}(\gamma[0,t])=t$.
Loewner's theorem states that in this case the maps $g_{t}$
satisfy the following differential equation
\begin{equation}\label{eq1}
  \partial_{t}{g_{t}(z)}=\frac{2}{g_{t}(z)-U(t)},~~g_{0}(z)=z,
\end{equation}
where $U(t)=g_{t}(\gamma(t))$, which is called the driving function for $\gamma$.

Conversely, consider a function $U(t):=\sqrt{\kappa}B(t)$, where $\kappa>0$ and $B(t)$
is a standard one-dimensional Brownian motion on $\mathbb{R}$
starting from $B(0)=0$. Then for
any $z\in \overline{\mathbb{H}}\setminus\{0\}$, the solution to (\ref{eq1}) exists
as long as $g_t(z)-U(t)$ stays away from $0$. Let $\tau(z)$ denote
the first time $\tau$ such that $\lim_{t\uparrow \tau}(g_t(z)-U(t))=0;  \tau(z)=\infty$ if
this never happens. Set $H_t:=\{z\in\mathbb{H}:\tau(z)>t\}$.
It is clear that $H_t$ is the set of points in $\mathbb{H}$ for which $g_t(z)$
is well-defined. Moreover, it is easy to verify that for
each $t\geq 0$, $g_t$ is a conformal map of $H_t$ onto $\mathbb{H}$, which
satisfies the hydrodynamic normalization (\ref{eq01}).

\begin{definit}The family of conformal maps $\{g_t:t\geq 0\}$ defined through (\ref{eq1}) is called the \textbf{ chordal $\mbox{SLE}_\kappa$
in $\mathbb{H}$}. The function $U(t)$ is called the driving function for the $\mbox{SLE}_\kappa$ process $\{g_t:t\geq 0\}$.
\end{definit}

Let $f(t,\cdot)$ denote the inverse of $g_t$, i.e., $f(t,\cdot):=g_t^{-1}$. The  trace $\gamma$ of $\mbox{SLE}_\kappa$
is defined by
\[\gamma(t):=\lim\limits_{z\rightarrow 0}f(t,z+U(t)),\]
where $z$ tends to $0$ within the upper half-plane $\mathbb{H}$. It was proved that $\gamma$ is a continuous path
in $\mathbb{H}$ from $0$ to $\infty$ (see
(\cite{RS05}($\kappa\neq8$) and \cite{LSW04}($\kappa=8$)). Moreover, it is easy to see that  $f(t,\cdot)$ satisfies the partial
differential equation
\begin{equation}\label{eq02}
\partial_tf(t,z)=-\partial_z f(t,z)\frac{2}{z-U(t)},\quad f(0,z)=z,z\in\mathbb{H}.
\end{equation}

Suppose that $D\subsetneq\mathbb{C}$
is a simply connected domain with two prescribed
distinct boundary points $z_0$ and $z_e$. Then the
Riemann mapping theorem implies that there is a conformal map $\psi:D\rightarrow\mathbb{H}$
satisfying $\psi(z_0)=0, \psi(z_e)=\infty$.
Let $\phi_t$ be the solution of the Loewner equation (\ref{eq1}) with initial condition $\phi_0(z)=\psi(z)$.
Then the process $\{\phi_t:t\geq 0\}$ is called the \textbf{chordal $\mbox{SLE}_\kappa$ in $D$} from $z_0$ to
$z_e$ under the map $\psi$. The scaling property of $\mbox{SLE}_\kappa$ \cite[Proposition 2.1(i)]{RS05}
implies that the behavior of
the process $\{\phi_t:t\geq 0\}$ does not depend on the choice of $\psi$,
so we will simply call $\{\phi_t:t\geq 0\}$ the chordal $\mbox{SLE}_\kappa$ process in $D$ from $z_0$ to
$z_e$, without mentioning the conformal map that maps $D$ to $\mathbb{H}$.
It is clear that $\phi_t=g_t\circ \psi$, where
$g_t$ is the solution of (\ref{eq1}) with initial condition $g_0(z)=z$. If $\tilde{\gamma}_t$
is the trace of the process $\{g_t\}$, then
the trace of the process $\{\phi_t\}$ is $\psi^{-1}(\tilde{\gamma}_t)$. Set $\gamma_t=\psi^{-1}(\tilde{\gamma}_t)$,
which  describes a cluster of
random curves in $D$ starting with $z_0$ and ending at $z_e$.

\subsection{Harmonic explorer}\label{haex}

In the subsection we briefly introduce the definition of harmonic explorer and state some related results.
Further details may be found in \cite{SS05}. It is worthwhile to point out that we will work in a slightly
different setting from  \cite{SS05} where larger and larger grid domains were considered,
while here we deal with rescaled grid domains.

As in Section \ref{intro}, let $TG^\epsilon$ denote the  triangular grid of the plane $\mathbb{C}$ with the lattice $V^\epsilon=\{\epsilon(m+ne^{i\pi/3}):m,n\in\mathbb{Z}\}$ for each positive number $\epsilon>0$.
Then the following result is well known as the Dirichlet problem of discrete harmonic function on $V^\epsilon$ (see \cite{SS05}).

\begin{lemma}\label{lem1}
Let $\Omega^\epsilon$ be any triangulation of a connected subregion of $TG^\epsilon$. If
$\hat{h}:V^\epsilon\cap\partial\Omega^\epsilon\rightarrow\mathbb{R}$ is a bounded function, then there is a unique
bounded function $h:V^\epsilon\cap \Omega^\epsilon\rightarrow\mathbb{R}$ which agrees with $\hat{h}$ in
$V^\epsilon\cap\partial\Omega^\epsilon$ and is harmonic at every vertex in $V^\epsilon\cap
(\Omega^\epsilon\setminus\partial\Omega^\epsilon)$.
\end{lemma}

The function $h$ in Lemma \ref{lem1} is said to be the discrete harmonic extension of $\hat{h}$.

Given a bounded Jordan domain $D\subset\mathbb{C}$ with two prescribed boundary points $a$ and $b$,
we define the $TG^\epsilon$ domain approximation of $D$, denoted by $D^\epsilon$, to be the largest
connected component of $TG^\epsilon$ contained in $D$. Here, $D^\epsilon$ is considered as
both a triangulation and a domain. Let $\hat{v}_0$ and $\hat{v}_e$ denote
the midpoints of $[v_0,u_0]$ and $[v_e,v_e]$ in $D^\epsilon$, respectively, which are closest to $a$ and $b$.
It is clear that $D^\epsilon$ tends to $D$ and that $\hat{v}_0$ and $\hat{v}_e$ also tend to $a$
and $b$ respectively, as $\epsilon\rightarrow 0$.

Write $\partial_+ D^\epsilon$ (resp.
$\partial_- D^\epsilon$) for the counterclockwise (resp. clockwise)
arc of $\partial D^\epsilon$ from $\hat{v}_0$ to $\hat{v}_e$.
Suppose further that $V^\epsilon\cap \partial_+ D^\epsilon$ is colored white,
and $V^\epsilon\cap \partial_- D^\epsilon$ colored  black.
Set $V_0^\epsilon:=V^\epsilon\cap \partial D^\epsilon$. We define $\hat{h}_0:V_0^\epsilon\rightarrow\{0,1\}$ to be $1$
on $V^\epsilon\cap\partial_+ D^\epsilon$, and $0$ on $V^\epsilon\cap\partial_- D^\epsilon$. Then Lemma \ref{lem1} implies
that there exists a unique discrete harmonic extension of $\hat{h}_0$, denoted by $h_0$.
\textbf{The harmonic explorer}, depending on the
triple $(D^\epsilon,\hat{v}_0,\hat{v}_e)$, is a random simple path from $\hat{v}_0$ to $\hat{v}_e$ described as follows.
Let $Z_1,Z_2,\hdots$ be independent and identically distributed  random variables, uniformly distributed in
the interval $[0,1]$. Let $\triangle_1^\epsilon\subset D^\epsilon$ be the triangle of $TG^\epsilon$ whose boundary contains $\hat{v}_0$,
and let $v_1$ be the vertex of $\triangle_1^\epsilon$ that is not on the edge containing $\hat{v}_0$.
Set $V_1^\epsilon:=V_0^\epsilon\cup\{v_1\}$. Let $v_1^r$ be the middle of the edge of $\triangle_1^\epsilon$ which
is on the counterclockwise arc from $\hat{v}_0$ to $v_1$, and let $v_1^l$ be the middle of the
edge of $\triangle_1^\epsilon$ which is on the clockwise arc from $\hat{v}_0$ to $v_1$.
If $Z_1\leq h_0(v_1)$, we let $\hat{v}_1:=v_1^l$, and otherwise $\hat{v}_1:=v_1^r$.
The beginning of a harmonic explorer is chosen as
the union of the two line segments from $\hat{v}_0$ to the center of $\triangle_1^\epsilon$ and then to $\hat{v}_1$.
Now a discrete function $\hat{h}_1:V_1^\epsilon\rightarrow \{0,1\}$ is defined by

\begin{equation*} \hat{h}_1(v):=
\begin{cases} \hat{h}_0(v),& \text{if $v\in V_0^\epsilon$},\\
\mathbf{1}_{Z_1\leq h_0(v_1)}, &\text{otherwise}.
\end{cases}
\end{equation*}
Thus we have defined the first step of harmonic explorer process.

Suppose that $n\geq 1$ and $\hat{v}_n\notin\partial D^\epsilon$. Then again by Lemma \ref{lem1}, there is a unique
discrete harmonic extension corresponding to $\hat{h}_n$, denoted by $h_n$. Let $\triangle_{n+1}^\epsilon\subset D^\epsilon$ be
the triangle of $TG^\epsilon$ containing $\hat{v}_n$ but not $\hat{v}_{n-1}$. Let $v_{n+1}$ be the
vertex of the edge of $\triangle_{n+1}^\epsilon$ which is not on the edge containing $\hat{v}_n$, and
let $V_{n+1}^\epsilon:=V_n^\epsilon\cup\{v_{n+1}\}$.
Let $\hat{v}_{n+1}^r$ and $\hat{v}_{n+1}^l$
be the two midpoints of edges of $\triangle_{n+1}^\epsilon$ containing $v_{n+1}$ which lie on the counterclockwise arcs of
$\partial \triangle_{n+1}^\epsilon$ from $\hat{v}_n$ to $v_{n+1}$ and from $v_{n+1}$ to $\hat{v}_n$, respectively.
If $Z_{n+1}\leq h_n(v_{n+1})$ let $\hat{v}_{n+1}:=\hat{v}_{n+1}^l$ and otherwise $\hat{v}_{n+1}:=\hat{v}_{n+1}^r$.
Then the next step of the harmonic explorer consists of segments from $\hat{v}_n$ to the center of $\triangle_{n+1}^\epsilon$ and
from the center of $\triangle_{n+1}^\epsilon$ to $\hat{v}_{n+1}$. Also, we define a discrete function $\hat{h}_n:V_{n+1}^\epsilon\rightarrow\{0,1\}$
by
\begin{equation*} \hat{h}_{n+1}(v):=
\begin{cases} \hat{h}_n(v),& \text{if $v\in V_n^\epsilon$},\\
\mathbf{1}_{Z_{n+1}\leq h_n(v_{n+1})}, &\text{otherwise}.
\end{cases}
\end{equation*}
It is easy to see that this procedure a.s. terminates  when $\hat{v}_n=\hat{v}_e$, and the harmonic explorer
constructed as above is a simple path from $\hat{v}_0$ to $\hat{v}_e$. Let $N$ denote the termination time, i.e.,
$\hat{v}_N=\hat{v}_e$. Then we have the following result. One can refer to \cite[Lemma 3.1]{SS05} for a proof.

\begin{lemma}\label{lem2} Let $h_n$ be defined as above. Then $h_n(v)$ is a martingale and $h_N(v)\in\{0,1\}$
for each $v\in D^\epsilon\cap V^\epsilon$.
\end{lemma}

Let $\gamma^\epsilon:[0,N]\rightarrow D^\epsilon\cup\{\hat{v}_0,\hat{v}_e\}$ be the harmonic explorer
path defined as above with the parameterization proportional to arc-length, where $\gamma^\epsilon(n)=\hat{v}_n$
for $n=0,1,\dots,N$. Let $\phi: D^\epsilon\rightarrow\mathbb{H}$ be a conformal map onto $\mathbb{H}$ with
$\phi(\hat{v}_0)=0$ and $\phi(\hat{v}_e)=\infty$. Then $\phi$ is unique up to a positive scaling,
and $\phi(\partial_- D^\epsilon)=(-\infty,0)$ and $\phi(\partial_+ D^\epsilon)=(0,\infty)$. Let $\tilde{\gamma}^\epsilon$ be the
path $\phi(\gamma^\epsilon)$, parameterized by the half plane capacity from $\infty$ in $\mathbb{H}$.

For convenience, we define a metric $\rho(\cdot,\cdot)$ on
$\overline{\mathbb{H}}$ by
$\rho(z,w)=|\varphi(z)-\varphi(w)|$, where $\varphi(z)=(z-i)/(z+i)$ is a conformal map from $\overline{\mathbb{H}}$ onto
the closed unit disk $\overline{\mathbb{U}}$. Then it is easy to see that for a given compact subset $K\subset\overline{\mathbb{H}}$,
there must be a constant $c>0$ depending only on $K$ such that
\begin{equation}\label{ineq1}
c^{-1}|z-z'|\leq \rho(z,z')\leq c|z-z'|
\end{equation}
for any $z, z'\in K$.
More details can be found in \cite[Section 3.2]{SS05}.

\subsection{Tip structure modulus}
In \cite{joh}  Viklund introduced the tip structure modulus for a radial Loewner curve, which is analogous
to Warschawski's boundary structure modulus \cite{wa}. In a similar manner, we will define the tip structure modulus
for a chordal Loewner curve.

Let $\gamma:[0,T]\rightarrow\overline{\mathbb{H}}$ be a curve with $\gamma(0)=0,\infty\notin\gamma[0,T]$, and for
$t\in[0,T]$, let $H_t$ be the unbounded connected component of $\mathbb{H}\setminus\gamma[0,t]$.
Then the curve $\gamma$ in $\mathbb{H}$, parameterized
by the half-plane capacity, is called an $\mathbb{H}$-Loewner curve if the following continuity condition holds:
for each $\epsilon>0$ there exists $\delta>0$ such that for all $s,t$ with $0<t-s<\delta$ there  is a crosscut
$\Gamma$ with $\mbox{diam}(\Gamma)<\epsilon$ that
separates $K_t\setminus K_s$ from $\infty$, where $K_t=\overline{\mathbb{H}\setminus H_t}$.
For each $t\geq 0$, the solution $f(t,\cdot)$ of the equation (\ref{eq02}),
where $U(t):[0,\infty)\rightarrow\mathbb{R}$ is a continuous function, is a conformal map from $\mathbb{H}$ onto some
simply connected domain $\tilde{H}_t\subset\mathbb{H}$.
\begin{definit}
The family $(f(t,z))_{t\geq 0}$ of conformal maps is said to be a Loewner chain.
A Loewner pair $(f,U)$ consists of a function $f(t,z)$ and a continuous function $U(t)$ such that $f$ is the solution to
the Loewner equation (\ref{eq02}) with $U(t)$ as the driving term.
\end{definit}

If $U$ is H\"{o}lder-$(1/2+\sigma)$ for some $\sigma>0$, then there
exists a curve $\gamma(t)$ such that $\tilde{H}_t$ is the unbounded connected component of $\mathbb{H}\setminus\gamma(t)$, and in this case we say
that the Loewner chain is generated by the Loewner curve $\gamma(t)$. Conversely,  given a Loewner curve $\gamma$, one can associate via
(\ref{eq02}) a unique driving term such that the Loewner chain $(f(t,\cdot))_{t\geq 0}$ in the Loewner pair $(f,U)$ is generated by $\gamma$.
In fact, the driving term is the preimage in $\partial\mathbb{H}$ under $f(t,\cdot)$ of the tip of the growing curve.

For a given domain $D^\epsilon$ defined in Section \ref{haex}, let $\phi:D^\epsilon\rightarrow\mathbb{H}$ be a conformal map onto $\mathbb{H}$
with $\phi(\hat{z}_0)=0$ and $\phi(\hat{z}_e)=\infty$.
Consider a chordal Loewner curve $\gamma:[0,T]\rightarrow D^\epsilon$ from $\hat{v}_0$ to $\hat{v}_e$. The conformal image
of $\gamma$ in $\mathbb{H}$ under $\phi$ is a $\mathbb{H}$-Loewner curve. We write $D_t^\epsilon$ for the connected component
of $D^\epsilon\setminus \gamma[0,t]$ containing $\hat{v}_e$.  For any $s,t\in[0,T]$, let $\gamma_{s,t}$
denote the curve determined by $\gamma(r),r\in[s,t]$. For each crosscut $\Gamma$ of $D_t^\epsilon$, denote by $J_\Gamma$ the
component of $D_t^\epsilon\setminus\Gamma$ with smaller diameter. For any $0\leq t\leq T$ and $\delta>0$, let $\mathcal{C}_{t,\delta}$ denote
the collection of crosscuts of $D_t^\epsilon$ with diameter at most $\delta$ which separate $\gamma(t)$ from $\hat{v}_e$ in $D_t^\epsilon$. For a crosscut
$\Gamma\in\mathcal{C}_{t,\delta}$, set
$s_\Gamma=\inf\{s>0:\gamma[t-s,t]\cap \bar{\Gamma}\neq\emptyset\}$, and $s_\Gamma=t$ if $\gamma$ never intersects $\bar{\Gamma}$.  Write $\gamma_\Gamma=\{\gamma(r):r\in[t-s_\Gamma,t]\}$.

\begin{definit}\label{def3}\quad For any $\delta>0$, the tip structure modulus of $(\gamma(t),t\in[0,T])$ in $D^\epsilon$, denoted
by $\eta(\delta)$, is defined by
\[ \eta(\delta)=\max\{\delta,\sup\limits_{t\in[0,T]}\sup\limits_{\Gamma\in \mathcal{C}_{t,\delta}}\mbox{diam}(\gamma_\Gamma)\}.\]
\end{definit}

Given $0<\delta\leq R$, (i) we shall say that $\gamma$ has a $(\delta,R)$-bottleneck in $D^\epsilon$ if there is $t\in[0,T]$ and
$w\in\partial D_t^\epsilon$ such that $\gamma(t)$ and $w$ can be connected by a crosscut $\Gamma_t$ of $D_t^\epsilon$ and $\mbox{diam}(J_{\Gamma_t})\geq
R$ while $\mbox{diam}(\Gamma_t)\leq \delta$; (ii) we will say that $\gamma$ has a nested $(\delta,R)$-bottleneck in $D^\epsilon$ if there is
a $t\in[0,T]$ and $\Gamma\in \mathcal{C}_{t,\delta}$ with $\mbox{diam}(\gamma_\Gamma)\geq R$.
Obviously, $\gamma(t),t\in[0,T]$ has no nested $(\delta,R)$-bottleneck in $D^\epsilon$ if and only if $\eta(\delta)\leq R$. This also
means that the curve $\gamma$ cannot visit a point $\zeta$, wander away, and then return to a point which is very near $\zeta$.

In the remaining text, if we consider the tip structure modulus for an $\mathbb{H}$-Loewner curve $\gamma(t)$, the distance
involved in Definition \ref{def3} refers to the metric $\rho$.

\section{Convergence rate for martingale observable}\label{rcmo}

In this section we will estimate the convergence rate for the martingale observable associated with the harmonic explorer
(see Proposition \ref{pro3-1}), which may be viewed as a quantitative version of \cite[Lemma 4.2]{SS05}.
Based on the triangle inequality, we reduce the required estimation to the following two tasks. One is
to estimate the difference between a discrete harmonic function and the corresponding non-discrete
harmonic function with the same boundary values (see Lemma \ref{lem3-2}). The other is to estimate the difference between two harmonic functions
with similar boundary values (see Lemma \ref{lem3-4}), by a result of Warschawski \cite{wa} (see Lemma \ref{lem3-3}).

For each $0\leq j<N$, let $\hat{D}_j^\epsilon:= D^\epsilon\setminus \cup_{l=1}^j\triangle_l^\epsilon$,
where $D^\epsilon$ and $\triangle_l^\epsilon$
are defined in Section \ref{haex}. Then $\hat{D}_j^\epsilon$ is a sub-triangulation of $D^\epsilon$ with
$\partial \hat{D}_j^\epsilon\cap V^\epsilon=V_j^\epsilon$. We write $\partial_+\hat{D}_j^\epsilon$ (resp.
$\partial_-\hat{D}_j^\epsilon$) for the counter-clockwise (resp. clockwise) arc from $\hat{v}_j$ to $\hat{v}_e$
of $\partial \hat{D}_j^\epsilon$. Then it follows from the construction of $h_j$
in Section 2.2 that $h_j$ is a discrete harmonic function on $\hat{D}_j^\epsilon\cap V^\epsilon$ which
satisfies the boundary condition that $h_j(v)=1, v\in\partial_+\hat{D}_j^\epsilon\cap V^\epsilon$
and $h_j(v)=0, v\in\partial_-\hat{D}_j^\epsilon\cap V^\epsilon$.
On the other hand,  we define the slit domains
\begin{displaymath}
D_{j}^\epsilon=D^\epsilon\backslash\bigcup_{i=1}^{j}[\gamma^\epsilon(i-1),\gamma^\epsilon(i)]=D^\epsilon\backslash\gamma^\epsilon[0,j]
\end{displaymath}
for any $1\leq{j}<N$, where $\gamma^\epsilon:[0,N]\rightarrow D^\epsilon\cup\{ \hat{v}_{0},\hat{v}_e\}$,
defined  in Section 2.2, is a harmonic explorer path satisfying $\gamma^\epsilon(j)=\hat{v}_j$
for each $j\in\{0,1,\cdots,N\}$. Let $\phi_{j}:D_{j}^\epsilon\rightarrow{\mathbb{H}}$ be the conformal map
from $D_{j}^\epsilon$ onto $\mathbb{H}$ satisfying $\phi_j\circ\phi^{-1}(z)-z\rightarrow 0$ as $z\rightarrow\infty$,
where $\phi$ is defined in Section 2.2. Denote by $W^\epsilon=W^\epsilon(t)$ the Loewner driving process
for $\tilde{\gamma}^\epsilon=\phi\circ \gamma^\epsilon$.
Then we have the following proposition.

\begin{prop}\label{pro3-1} For every $j<N$ and for any compact subset $K\subset \hat{D}_j^\epsilon$, there exists
a constant $C$ depending only $K$ such that for each vertex $v\in K\cap V^\epsilon$,
\begin{equation}\label{eq11} |h_j(v)-\tilde{h}(\phi_j(v)-W^\epsilon(t_j))|<C\epsilon^{1/2}
\end{equation}
for $\epsilon$ sufficiently small, where $\tilde{h}(z)=1-(1/\pi)\arg z$.
\end{prop}

To prove Proposition \ref{pro3-1}, we need the following lemmas.

\begin{lemma}\label{lem3-2} For arbitrarily small $\epsilon>0$, let $Q^\epsilon\subset TG^\epsilon$ be a grid
bounded domain with mesh size $\epsilon$. Let $h$ be a discrete
harmonic function on $Q^\epsilon\cap V^\epsilon$, and let $\bar{h}$ be a harmonic function on a bounded domain containing $\overline{Q}^\epsilon$.
If $h(v)=\bar{h}(v)$ for each boundary vertex $\partial Q^\epsilon\cap V^\epsilon$,
then there exists a universal constant $C>0$ such that
\begin{equation}\label{equ3-1} |h(v)-\bar{h}(v)|<C\epsilon^2
\end{equation}
for any vertex $v\in Q^\epsilon\cap V^\epsilon$.
\end{lemma}

\noindent\textbf{Proof.} First, for $\bar{h}$ which is viewed as a discrete function on
$Q^\epsilon\cap V^\epsilon$ we will show
\begin{equation}\label{eq3-2}\Delta \bar{h}(v)=O(\epsilon^4)
\end{equation}
for each $v\in Q^\epsilon\cap V^\epsilon$, where $\Delta$ is the discrete Laplacian,
that is,
\[\Delta \bar{h}(v)=\frac{1}{6}\sum\limits_{k=0}^5[\bar{h}(v+\epsilon e^{i\frac{k\pi}{3}})-\bar{h}(v)].\]
Indeed, since $\bar{h}$ is harmonic, we can  expand $\bar{h}$ in power series about each $z_0\in Q^\epsilon$ to get
\[\bar{h}(z)=a_0+a_1(z-z_0)+a_2(z-z_0)^2+a_3(z-z_0)^3+O(|z-z_0|^4),\]
where $a_j(j=0,1,2,3)$ are constants associated with $z_0$. This implies
\[\Delta \bar{h}(v)=\frac{1}{6}\sum\limits_{k=0}^5(\tilde{a}_1\epsilon e^{i\frac{k\pi}{3}}+\tilde{a}_2\epsilon^2e^{i\frac{2k\pi}{3}}
+\tilde{a}_3\epsilon^3e^{ik\pi})+O(\epsilon^4),\]
where $\tilde{a}_j(j=0,1,2,3)$ are constants associated with $v$. A straightforward computation gives that
\[\sum\limits_{k=0}^5 e^{i\frac{jk\pi}{3}}=0\]
for $j=1,2,3$. Hence, we deduce that (\ref{eq3-2}) holds.

Next, consider the function $f(v)=h(v)-\bar{h}(v)+\beta|v/\epsilon|^2$, where $\beta\in(0,\epsilon^2)$ is a
function of $\epsilon$.  We intend to determine how large $\beta$ is so that $f$ has no maximum  in
$\text{int}(Q^\epsilon\cap V^\epsilon)$, the set of interior vertices of $Q^\epsilon$. Now suppose that $f$ has
a maximum at some vertex $v_0\in \text{int}(Q^\epsilon\cap V^\epsilon)$.
Then for $k=0,1,\dots,5$, we get from the definition of $f$ that
\begin{equation}\label{eq3-3}
h(v_0+\epsilon e^{i\frac{k\pi}{3}})\leq p_k,
\end{equation}
where
\begin{equation}\label{eq3-4}
p_k=h(v_0)+\bar{h}(v_0+\epsilon e^{i\frac{k\pi}{3}})-\bar{h}(v_0)-\beta|v_0/\epsilon+e^{i\frac{k\pi}{3}}|^2
+\beta|v_0/\epsilon|^2.
\end{equation}
It follows from (\ref{eq3-4}) and (\ref{eq3-2}) that
\[\sum\limits_{k=0}^5 p_k-6h(v_0)=-\beta\sum\limits_{k=0}^5(|v_0/\epsilon+e^{i\frac{k\pi}{3}}|^2
-|v_0/\epsilon|^2)+O(\epsilon^4).\]
A straightforward computation gives that
\[\sum\limits_{k=0}^5(|v_0/\epsilon+e^{i\frac{k\pi}{3}}|^2
-|v_0/\epsilon|^2)=6.\]
Thus we obtain
\begin{equation}\label{eq3-5}\sum\limits_{k=0}^5 p_k-6h(v_0)=-6\beta+O(\epsilon^4).
\end{equation}
Since $h$ is discrete harmonic, i.e., $h(v_0)=(1/6)\sum_{k=0}^5h(v_0+\epsilon e^{i\frac{k\pi}{3}})$,
we obtain from (\ref{eq3-5}) that
\[\sum\limits_{k=0}^5(p_k-h(v_0+\epsilon e^{i\frac{k\pi}{3}}))=-6\beta+O(\epsilon^4).\]
This, combined with (\ref{eq3-3}), implies
\[\beta=O(\epsilon^4).\]
It means that if we choose $\beta=C\epsilon^4$ with $C>0$ a sufficiently large constant, then $f$ will have no
maximum in $\text{int}(Q^\epsilon\cap V^\epsilon)$. In this case, we have $f(v)=\beta|v/\epsilon|^2=O(\epsilon^2)$
for $v\in \partial Q^\epsilon\cap V^\epsilon$.
Hence we get $f(v)\leq O(\epsilon^2)$ for all $v\in \text{int}(Q^\epsilon\cap V^\epsilon)$. This, combined with the definition of $f$, implies
\begin{equation}\label{eq3-6}h(v)-\bar{h}(v)\leq O(\epsilon^2)
\end{equation}
for any $v\in Q^\epsilon\cap V^\epsilon$.

Finally, we consider the function $g(v)=\bar{h}(v)-h(v)+\gamma|v/\epsilon|^2$, where $\gamma\in(0,\epsilon^2)$ is a
function of $\epsilon$. Following almost identical lines in the previous paragraph, we conclude that
\begin{equation}\label{eq3-7} \bar{h}(v)-h(v)\leq O(\epsilon^2)
\end{equation}
for each $v\in Q^\epsilon\cap V^\epsilon$. Thus it follows from (\ref{eq3-6}) and (\ref{eq3-7}) that (\ref{equ3-1}) holds.
$\hfill\square$

\begin{lemma}\label{lem3-3} Suppose that $D$ and $D'$ are two simply connected  bounded domains such
that $D'\subset D$ and $z=0$ lies in $D'$. Let $\varphi$ and $\psi$ be the two conformal maps from $D$ and $D'$,
respectively, onto the unit disk $\mathbb{U}$ normalized by $\varphi(0)=0,\varphi'(0)>0$ and $\psi(0)=0,
\psi'(0)>0$. If $K\subset D'$ is a compact subset and $\mbox{dist}(z,\partial D)\leq \epsilon$ for all
$z\in \partial D'$, then there exists a constant $C=C(K)>0$ depending only on $K$ such that
\begin{equation*}\label{eq3-8}
|\varphi(z)-\psi(z)|\leq C\epsilon^{1/2}
\end{equation*}
for any $z\in K$.
\end{lemma}

\noindent\textbf{Proof.} It follows from \cite[Lemmas 3 and 5]{wa}. $\hfill\square$

\begin{lemma}\label{lem3-4} Suppose that $D$ and $D'$ are two simply connected  bounded domains such
that $D'\subset D$ and $z=0$ lies in $D'$. Let $u$ and $\tilde{u}$ be the two harmonic functions on $D$ and $D'$,
respectively, with the following boundary values
\begin{equation*} u(z)=
\begin{cases} 1,& \text{if $z\in L$, where $L\subset\partial D$ is a sub-arc of $\partial D$ with  positive length },\\
0, &\text{if $z\in (\partial D\setminus L)$}
\end{cases}
\end{equation*}
and
\begin{equation*} \tilde{u}(z)=
\begin{cases} 1,& \text{if $z\in \tilde{L}$, where $\tilde{L}\subset\partial D'$ is a sub-arc of $\partial D'$ with positive length },\\
0, &\text{if $z\in (\partial D'\setminus \tilde{L})$}.
\end{cases}
\end{equation*}
If $K\subset D'$ is a compact subset
and $\mbox{dist}(z,\partial D)\leq \epsilon$ for all
$z\in \partial D'$, then there exists a constant $C=C(K)>0$ depending only on $K$ such that
\begin{equation}\label{eq3-9}
|u(z)-\tilde{u}(z)|\leq C\epsilon^{1/2}
\end{equation}
for any $z\in K$.
\end{lemma}

\noindent\textbf{Proof.} We first claim that there exist two conformal maps $\varphi:D\rightarrow S$
and $\psi:D'\rightarrow S$, where $S\subset\mathbb{C}$ is some simply connected domain,
such that the following results hold:
\begin{equation}\label{eq3-9-1}\mbox{Re}(\varphi(z))=u(z),\mbox{Re}(\psi(z))=\tilde{u}(z),
\end{equation}
and
\begin{equation}\label{eq3-9-2} \varphi(0)=\psi(0),\varphi'(0)>0,\psi'(0)>0.
\end{equation}
Indeed, since $u$ is a solution of the Laplace equation $u_{xx}+u_{yy}=0$ on $D$,
there exists a companion solution $v(z)\in D$ to the Cauchy-Riemann equation $u_x=v_y,u_y=-v_x$. Moreover,
$v(z)$ is unique up to an additive scalar. This implies that there is an analytic function $\varphi$ on $D$
such that $\varphi(z)=u(z)+iv(z)$. On the other hand, the assumption on the boundary value of $u$ gives that for any $z\in D$,
$u(z)$ is the harmonic measure of $L$ from the point $z$ inside $D$. This yields that $u(z_1)\neq u(z_2)$
for any $z_1\neq z_2\in D$. Hence we obtain that $\varphi$ is a univalent function on $D$.  The same arguments carried
to $\tilde{u}$ give that there exists a univalent function $\psi$ on $D'$
such that $\psi(z)=\tilde{u}(z)+i\tilde{v}(z)$, where $\tilde{v}(z)\in D'$ is
a companion solution to the corresponding Cauchy-Riemann equation.
Since a nonconstant analytic function maps a domain onto another one, both the images $\varphi(D)$ and $\psi(D')$
are simply connected domains. At the same time, notice that
$\{\text{Re}\varphi(z):z\in \partial D\}=\{\text{Re}\psi(z):z\in \partial D'\}=\{0,1\}$.  So we deduce that $\varphi(D)=\psi(D')
=\{w=u+iv\in\mathbb{C}:0<u<1,-\infty<v<\infty\}\triangleq S$. Both $\varphi$ and $\psi$ are
two conformal maps from $D$ and $D'$ onto $S$, respectively, such that the result (\ref{eq3-9-1}) holds.  By the Riemann mapping theorem
we can further choose $\varphi$ and $\psi$ satisfying
(\ref{eq3-9-2}), too. This implies that the claim holds.

Next, again applying the Riemann mapping theorem to $S$ gives that there is a unique conformal
map $f:S\rightarrow\mathbb{U}$ satisfying $f(\varphi(0))=0, f'(\varphi(0))>0$. Write $\Phi=f\circ\varphi$
and $\Psi=f\circ\psi$. Then $\Phi$ and $\Psi$ are  two conformal maps from $D$ and $D'$,
respectively, onto the unit disk $\mathbb{U}$ normalized by $\Phi(0)=0,\Phi'(0)>0$ and $\Psi(0)=0,
\Psi'(0)>0$.  Thus it follows from Lemma \ref{lem3-3} that there exists a constant $C>0$ depending
only on $K$ such that
\begin{equation}\label{eq3-a}|\Phi(z)-\Psi(z)|\leq C\epsilon^{1/2}
\end{equation}
for any $z\in K$.

Finally, note that an analytic function is an open map. So we get
that $\varphi(K)\cup\psi(K)$ is a compact subset of $S$. Thus applying Koebe's distortion theorem  to $f$,
combined with (\ref{eq3-a}), gives that
\begin{equation}\label{eq3-10} |\varphi(z)-\psi(z)|<C|f(\varphi(z))-f(\psi(z))|\leq C\epsilon^{1/2}
\end{equation}
for any $z\in K$, where $C>0$ is a constant depending only on $K$. Since $u(z)=\mbox{Re}(\varphi(z))$ and
$\tilde{u}(z)=\mbox{Re}(\psi(z))$, we conclude from (\ref{eq3-10}) that (\ref{eq3-9}) holds. $\hfill\square$

\

\noindent\textbf{Proof of Proposition \ref{pro3-1}.} Let  $\bar{h}_j$ be the harmonic function on $\hat{D}_j^\epsilon$
satisfying the boundary conditions: $\bar{h}_j(z)=1, z\in \partial_+D_j^\epsilon$ and $\bar{h}_j(z)=0,
z\in \partial D_j^\epsilon\setminus\partial_+D_j^\epsilon$.
Then by the triangle inequality we have
\begin{equation}\label{eq3-12} |h_j(v)-\tilde{h}(\phi_j(v)-W^\epsilon(t_j))|\leq |h_j(v)-\bar{h}_j(v)|
+|\bar{h}_j(v)-\tilde{h}(\phi_j(v)-W^\epsilon(t_j))|,
\end{equation}
We will estimate each summand on the right-hand side of the inequality in (\ref{eq3-12}).

For the first term, we notice that  $\bar{h}_j$ may be viewed as the real part of one analytic
function $F(z)$ on $\hat{D}_j^\epsilon$. It is clear that $F(z)$ can be extended analytically
to $\bar{F}(z)$ on a domain containing $\overline{\hat{D}_j^\epsilon}$ such that
$\bar{F}(z)=F(z)$ for $z\in \overline{\hat{D}_j^\epsilon}$. Setting $\bar{H}_j=\mbox{Re}\bar{F}_j$,
one has $\bar{h}_j=\bar{H}_j$ on $\overline{\hat{D}_j^\epsilon}$. So applying Lemma \ref{lem3-2} to $h_j$ and $\bar{H}_j$ gives
\begin{equation}\label{eq3-13}
|h_j(v)-\bar{h}_j(v)|=|h_j(v)-\bar{H}_j(v)|<C\epsilon^2,
\end{equation}
where $C>0$ is a universal constant.

Next, consider the second term  on the right-hand side of the inequality in (\ref{eq3-12}). Obviously, we can take a
fixed point $z_0\in\hat{D}_j$ for all $j<N$. Without loss of generality we may assume that
$z_0=0$. This, combined with the definitions of $\bar{h}_j$ and $\tilde{h}(\phi_j(\cdot)-W^\epsilon(t))$, implies that
the conditions of Lemma \ref{lem3-4} are satisfied when $u$ and $\tilde{u}$ are replaced by $\bar{h}_j$
on $\hat{D}_j^\epsilon$ and $\tilde{h}(\phi_j(\cdot)-W^\epsilon(t))$ on $D_j^\epsilon$,
respectively. Thus we apply Lemma \ref{lem3-4} to $\bar{h}_j$ and $\tilde{h}(\phi_j(\cdot)-W^\epsilon(t))$ to obtain that
\begin{equation}\label{eq3-14}
|\bar{h}_j(v)-\tilde{h}(\phi_j(v)-W^\epsilon(t_j))|<C\epsilon^{1/2}
\end{equation}
for each vertex $v\in K$, where $C>0$ is a constant depending only on $K$.

Therefore, from (\ref{eq3-12}),(\ref{eq3-13}) and (\ref{eq3-14}) we conclude that (\ref{eq11}) holds. $\hfill\square$

\section{Moment estimates for increments of driving function} \label{meidf}

In this section, based on Proposition \ref{pro3-1} we will derive the moment estimates for increments of the driving function
of harmonic explorer process. The idea of proof is similar to  \cite[proposition 4.1]{SS05}.
However, the main difference is that we provide a decaying rate, which may be viewed as a quantitative
version of \cite[proposition 4.1]{SS05}. In addition, the related sets are different.
Instead of larger and larger grid domains discussed in \cite{SS05}, we work on compact
subsets of domain $D$ since we consider rescaled grid domains, i.e., the mesh of grid tends to zero.

As in Section \ref{haex}, let $\tilde{\gamma}^\epsilon$ denote the image of the harmonic explorer $\gamma^\epsilon$
under the conformal map $\phi$, parameterized by the half plane capacity from $\infty$
in $\mathbb{H}$, and $W^\epsilon=W^\epsilon(t)$ the Loewner driving process for $\tilde{\gamma}^\epsilon$.
Let $t_n$ be the half plane capacity of the curve $\tilde{\gamma}^\epsilon[0,n])$ from $\infty$ in $\mathbb{H}$
for each $n<N$. Then we have the following moment estimates for increments of $W^\epsilon$.

\begin{prop}\label{pro4-1} Fix any integer $n\geq 0$. On the event $E_1=E_1(n):=\{n<N\}$, define
\[m:=\min{\{k>n:(t_{k}-t_{n})\vee(W^\epsilon(t_k)-W^\epsilon(t_n))^{2}\geq{\epsilon^{1/3}}\}}.\]
Write $p_n=(\phi_n)^{-1}(i+W^\epsilon(t_n))$ and let $E_2=E_2(n)$ denote the event that
$K\subset \hat{D}_n^\epsilon$ is a compact subset containing $p_n$, where $\phi_n$ and $\hat{D}_n^\epsilon$
are defined  in Section \ref{rcmo}. Then
for sufficiently small $\epsilon$, it holds that
\begin{equation}\label{eq2}
\mid\mathbb{E}[W^\epsilon(t_m)-W^\epsilon(t_{n})\mid{\gamma^\epsilon[0,n]}]\mid\leq C\epsilon^{1/2}
\end{equation}
and
\begin{equation}\label{eq3}
\mid\mathbb{E}[(W^\epsilon(t_m)-W^\epsilon(t_{n}))^{2}-4(t_m-t_{n})\mid{\gamma^\epsilon[0,n]}]\mid\leq C \epsilon^{1/2}
\end{equation}
on the event $E_1\cap E_2$, where $C>0$ is a constant depending only on $K$.
\end{prop}

\noindent{\bf Proof.} Start with $E_1$.  We claim that if $K\subset \hat{D}_n^\epsilon$ is a compact subset
of $\hat{D}_n^\epsilon$,
then $K$ is also a compact subset of $\hat{D}_m^\epsilon$ for $\epsilon$ sufficiently small.
Indeed, from the compactness of $K$ it is easy to see  that $d:=\mbox{dist}(K,\partial \hat{D}_n^\epsilon)>0$.
On the other hand, the definition of $m$ gives  that $(t_{m-1}-t_n)^{1/2}\leq \epsilon^{1/6}$ and $|W(t_{m-1})-W(t_{n})|\leq \epsilon^{1/6}$.
It follows from \cite[Lemma 2.1]{LSW04} that $\mbox{diam}(\gamma^\epsilon_\phi[n,m-1])\leq c\epsilon^{1/6}$, where
$c>0$ is a universal constant.
Thus applying Koebe's distortion theorem to $\phi^{-1}$ gives that
there exists an $0<\epsilon_0<d/4$ such that
\[\mbox{diam}(\gamma^\epsilon[n,m-1])<d/4\]
for all $\epsilon<\epsilon_0$. This, combined with the definition of $\gamma^\epsilon(t_j)$,
implies that $\mbox{diam}(\gamma^\epsilon[n,m])\leq\mbox{diam}(\gamma^\epsilon[n,m-1])+\mbox{diam}(\gamma^\epsilon[m-1,m])<d/4+\epsilon<d/2$.
Hence we obtain that $\mbox{dist}(K,
\partial \hat{D}_m^\epsilon)\geq d/2>0$, which implies that the claim holds.

Next, assume $E_2$ occurs, and fix some vertex $v_0\in K\cap V^\epsilon$. The above claim  gives
that we can apply (\ref{eq11}) with $j=n,m$ and $v=v_0$.
Moreover, Lemma \ref{lem2} gives that $h_j(v)$ is a martingale
for any $v\in D^\epsilon\cap V^\epsilon$.  So we get
\[\mathbb{E}[h_m(\upsilon_0)|\gamma^\epsilon[0,n]]=h_n(\upsilon_0).\]
This, combined with (\ref{eq11}),  implies
\begin{equation}\label{eq3-1}\mathbb{E}[\widetilde{h}(\phi_{m}(\upsilon_{0})-W^\epsilon(t_{m}))\mid{\gamma^\epsilon[0,n]}]
=\widetilde{h}(\phi_{n}(\upsilon_{0})-W^\epsilon(t_{n}))+O{(\epsilon^{1/2})}.
\end{equation}
In addition, from (4.5) in \cite{SS05} we obtain that
\begin{equation}\label{eq6}
\mid{W^\epsilon(t)-W^\epsilon(t_{n})}\mid=O(\epsilon^{1/6}),
\quad\mid{t_{m}-t_{n}}\mid=O(\epsilon^{1/3})
\end{equation}
for any $t\in[t_{n},t_{m}]$.

Set $z_{t}:=g_{t}\circ\phi{(\upsilon_{0})}$, where $g_t$ is the solution to equation (\ref{eq1}) with $U(t)=W^\epsilon(t)$.
Then we have
$\phi_{j}(\upsilon_{0})=z_{t_{j}}$, where $\phi_{j}$ is defined in Section \ref{rcmo}.
This yields that $z_{t_m}$ can be obtained from $z_{t_n}$ according to (\ref{eq1}) with $t$
between $t_n$ and $t_m$. Note that $\mbox{Im} z_t$ has a constant positive lower bound for $t\in [t_n,t_m]$.
So from (\ref{eq6}) we can deduce that
$$\frac{2}{z_{t}-W^\epsilon(t)}=\frac{2}{z_{t_{n}}-W^\epsilon(t_{n})}+O(\epsilon^{1/6})$$
for all $t\in [t_n,t_m]$.
Integrating the above equality over $[t_{n},t_{m}]$,
together with (\ref{eq1}), implies
\begin{equation}\label{equ6-1}z_{t_{m}}-z_{t_{n}}=\phi_{m}(\upsilon_{0})-\phi_{n}(\upsilon_{0})
=\frac{2(t_{m}-t_{n})}{\phi_{n}(\upsilon_{0})-W^\epsilon(t_{n})}+O{(\epsilon^{1/2})}.
\end{equation}

Finally, we let
\[G(z,W^\epsilon):=\widetilde{h}{(z-W^\epsilon)}=1-(1/\pi)\arg(z-W^\epsilon)=1-\frac{1}{2\pi i}\log\frac{z-W^\epsilon}{\bar{z}-W^\epsilon}.\]
Then applying a Taylor expansion to
$G$ at $(z_{t_{n}}, W^\epsilon(t_{n}))$ gives
\begin{eqnarray*}\widetilde{h}(\phi_{m}(\upsilon_{0})&-&W^\epsilon(t_{m}))-\widetilde{h}(\phi_{n}(\upsilon_{0})-W^\epsilon(t_{n}))\\
&=&\partial_{z}{G_{(z_{t_{n}}, W^\epsilon(t_{n}))}{(z_{t_{m}}-z_{t_{n}})}}+\partial_{W^\epsilon}{G_{(z_{t_{n}}, W^\epsilon(t_{n}))}{(W^\epsilon(t_{m})-W^\epsilon(t_{n}))}}
\\&&+(1/2)\partial_{W^\epsilon}^{2}{G_{(z_{t_{n}}, W^\epsilon(t_{n}))}{(W^\epsilon(t_{m})-W^\epsilon(t_{n}))^{2}}}+O{(\epsilon^{1/2})}.
\end{eqnarray*}
Taking the conditional expectations to the both sides of the above equality given $\gamma^\epsilon[0,n]$,
combined with (\ref{eq3-1}) and (\ref{equ6-1}), implies that
\begin{eqnarray}\label{eq8}
O{(\epsilon^{1/2})}&=&\mbox{Im}((\phi_{n}(\upsilon_{0})-W^\epsilon(t_{n}))^{-2})\mathbb{E}[4(t_{m}-t_{n})
-(W^\epsilon(t_{m})-W^\epsilon(t_{n}))^{2}\mid{\gamma^\epsilon[0,n]}]
\nonumber\\&&-2\mbox{Im}((\phi_{n}(\upsilon_{0})-W^\epsilon(t_{n}))^{-1})
\mathbb{E}[(W^\epsilon(t_{m})-W^\epsilon(t_{n}))\mid{\gamma^\epsilon[0,n]}].
\end{eqnarray}

Let $v_1\in \hat{D}_n^\epsilon\cap V^\epsilon$ be a vertex closest to $p_n$,
and let $v_2\in \hat{D}_n^\epsilon\cap V^\epsilon$ be a vertex closest to
$\phi_n^{-1}(i+W^\epsilon(t_n)+1/100)$. The Koebe distortion theorem
gives $|\phi_n(v_1)-i-W^\epsilon(t_n)|=|\phi_n(v_1)-\phi_n(p_n)|=O(\epsilon)$ and $|\phi_n(v_2)-i-W^\epsilon(t_n)-1/100|=O(\epsilon)$.
This, combined with the fact that $\mbox{Im} \phi_n(p_n)=1$, implies that both
$\text{Im}\phi_n(v_1)$ and $\text{Im}\phi_n(v_2)$ have a positive lower
bound for $\epsilon$ sufficiently small. Hence we get that $\mbox{Im}(\phi_n(v_1)-W^\epsilon(t_n))^{-1}
\neq \mbox{Im}(\phi_n(v_2)-W^\epsilon(t_n))^{-1}$ and $\mbox{Im}(\phi_n(v_1)-W^\epsilon(t_n))^{-2}
\neq \mbox{Im}(\phi_n(v_2)-W^\epsilon(t_n))^{-2}$, each of which is bounded away from zero.
On the other hand, since $p_n\in K$, we have $v_1,v_2\in K$ for $\epsilon$ sufficiently small.
Thus we can apply (\ref{eq8}) with $v_0$ replaced by each of $v_1,v_2$,
which produces two linearly independent equations of the variables $\mathbb{E}[4(t_{m}-t_{n})
-(W^\epsilon(t_{m})-W^\epsilon(t_{n}))^{2}\mid{\gamma^\epsilon[0,n]}]$ and
$\mathbb{E}[(W^\epsilon(t_{m})-W^\epsilon(t_{n}))\mid{\gamma^\epsilon[0,n]}]$. Solving these two linear equations  implies
(\ref{eq2}) and (\ref{eq3}), where $C>0$ is a constant depending only on $K$.
This completes the proof. $\hfill\square$

\section{Convergence rate for driving function} \label{ecrdf}

In this section, using Proposition \ref{pro4-1} and the Skorokhod embedding theorem we will prove Theorem \ref{thm1}, i.e.,
derive the rate of convergence of the driving function for harmonic explorer to the Brownian motion $B(4t)$, which is the Loewner
driving function for the chordal $\mbox{SLE}_4$. Although the structure of the proof is similar to that of \cite[Theorem 1.1]{BVK},
there are the following differences. First, the driving functions investigated are different. The driving function
discussed in \cite{BVK} is induced by the loop-erased random walk; while
the driving function considered here is produced by the harmonic explorer process. So a lot of details involved in our proof are different
from those in \cite{BVK}, in particular to ensure that the conditions of Proposition \ref{pro4-1} are satisfied.
Secondly, the results obtained are different. The exponent of convergence rate given in \cite{BVK} is $1/24$; whereas the exponent of convergence
rate we here obtain is $1/12$. This shows that the convergence rate of the latter is faster than one of the
former.

In order to prove Theorem \ref{thm1}, we need the following lemmas. First, the Skorokhod embedding theorem for martingale
is needed; see \cite[Theorem A.1]{HH} or \cite{ry} for a proof.

\begin{lemma}\label{lem5-2}(Skorokhod embedding)  Assume that $(M_j)_{j\leq J}$ is an $(\mathcal{F}_{j})_{j\leq J}$ martingale, with $\|
M_j-M_{j-1}\|_{\infty}\leq \delta$ and $M_0=0$ a.s. Then there are stopping times $0=\tau_0\leq\tau_1\leq
\cdots,\leq \tau_J$ for a standard Brownian motion $(B(t),t\geq 0)$, such that $(M_0,M_1,\dots,M_J)$ and $(B(\tau_0),
B(\tau_1),\dots,B(\tau_J))$ have the same law. Moreover, it holds that for $j=0,1,\dots, J-1$,
\begin{equation}\label{eq9-1}
\mathbb{E}[\tau_{j+1}-\tau_j|B[0,\tau_j]]=\mathbb{E}[(B(\tau_{j+1})-B(\tau_j))^2|B[0,\tau_j]],
\end{equation}
\begin{equation}\label{eq9-2}
\mathbb{E}[(\tau_{j+1}-\tau_j)^q|B[0,\tau_j]]\leq C_q\mathbb{E}[(B(\tau_{j+1})-B(\tau_j))^{2q}|B[0,\tau_j]]
\end{equation}
where $C_q>0$ is a constant, and
\begin{equation}\label{eq9-3}
\tau_{j+1}\leq\inf\{t\geq\tau_j:|B(t)-B(\tau_k)|\geq\delta\}.
\end{equation}
\end{lemma}

Next, we will use the following inequality for martingale difference sequence; the proof can be found
in \cite[Lemma 1]{HA}.

\begin{lemma}\label{lem5-3} Suppose that $X_j,j=1,2,\dots,J$, is a martingale difference sequence with respect
to the filtration $\mathcal{F}_j$. Then for any $a,b,c>0$,
\begin{eqnarray*}
\mathbb{P}\Big(\max\limits_{1\leq j\leq J}|\sum_{k=1}^jX_k|\geq a\Big)&\leq&\sum\limits_{j=1}^J\mathbb{P}(|X_j|>a)
+2\mathbb{P}\Big(\sum\limits_{j=1}^J\mathbb{E}[X_j^2|\mathcal{F}_{j-1}|]>b\Big)\\
&&+2\exp\{ab^{-1}(1-\log(abc^{-1}))\}.
\end{eqnarray*}
\end{lemma}

Meanwhile, we will also employ the following result about the modulus of continuity of Brownian motion;
see \cite[Lemma 1.2.1]{CR} for the proof.

\begin{lemma}\label{lem5-4} Suppose that $B(t),t\geq 0$, is a standard Brownian motion. Then for any $\delta>0$ there
exists a constant $C=C(\delta)>0$ such that
\[\mathbb{P}\Big(\sup\limits_{t\in[0,T-\lambda]}\sup\limits_{s\in[0,\lambda]}|B(t+s)-B(t)|\leq \mu\sqrt{\lambda}\Big)
\geq 1-\frac{CT}{\lambda}e^{-\frac{\mu^2}{2+\delta}}\]
for any positive numbers $\mu,T$ and $0<\lambda<T$.
\end{lemma}

\noindent\textbf{Proof of Theorem \ref{thm1}.}  For any fixed $R>0$, we define $\tilde{T}\leq T$ by $\tilde{T}:=(1/24)\sup\{t\in[0,24T]:|W^\epsilon(t)|\leq
R\}$. Set $I:=\{n\in\mathbb{N}:t_n\leq 24\tilde{T}\}$ for $\epsilon$ sufficiently small.
In order to apply Proposition \ref{pro4-1}, we need to verify
$E_1\cap E_2$ for each $n\in I$. It follows from \cite[Lemma 2.1]{LSW04} that $t_N=\infty$ or $\{W(t):t \in[0,t_N]\}$
is unbounded, which implies $E_1$ for any $n\in I$.
Note that $\phi_n(p_n)=i+W^\epsilon(t_n)$ and $g_{t_n}=\phi_n\circ\phi^{-1}$, so we get $g_{t_n}\circ\phi(p_n)=i+W^\epsilon(t_n)$.
We claim that there exists a compact subset $\tilde{K}\subset\mathbb{H}$,  depending only on $R$ and $\tilde{T}$,
such that $\phi(p_n)\in \tilde{K}$ holds for any $n\in I$. In fact, $g_t\circ \phi(p_n)$ satisfies
(\ref{eq1}) starting from $\phi(p_n)$ at $t=0$ to $i+W(t_n)$ at $t=t_n$. The monotonicity of $\mbox{Im} g_t$ with respect
to $t$ gives that $\mbox{Im} g_t\circ\phi(p_n)\geq 1$ for each $t\in [0,t_n]$. This, combined with (\ref{eq1}), implies
that $|\partial_tg_t\circ\phi(p_n)|=O(1)$. Hence we deduce that $|\phi(p_n)|\leq 1+|W^\epsilon(t)|+O(\tilde{T})\leq 1+R+O(\tilde{T})$.
Let $\tilde{K}:=\{z\in\mathbb{C}:\mbox{Im}(z)\geq 1, |z|\leq O(\tilde{T}+R)\}$. Then it is clear that $\tilde{K}$
is a compact subset of $\mathbb{H}$
containing $\phi(p_n)$ for all $n\in I$, which implies that the claim holds. Let $K=\phi^{-1}(\tilde{K})$. Then $K$ is
a compact subset of $D_n\subset D$ for all $n\in I$, since $\phi^{-1}$ is an open map.
Moveover, $K$ contains all $p_n$ for $n\in I$. It is easy to see that for each $n\in I$, $K$ is also a compact subset of $\hat{D}_n^\epsilon$
containing $p_n$ for $\epsilon$ sufficiently small. This implies that $E_2$ occurs for each $n\in I$. So the conditions
of Proposition \ref{pro4-1} are satisfied for any $n\in I$.

Set $m_0:=n=0$ and $m_1:=m$, where $m$ is defined in Proposition \ref{pro4-1}.
Inductively on $j=1,2,\dots$, define
\[m_{j+1}=\min\{k>m_j:|t_k-t_{m_j}|\geq \epsilon^{1/3}\vee |W^\epsilon(t_k)-W^\epsilon(t_{m_j})|\geq \epsilon^{1/6}\}.\]
Write $J=\lceil 12\tilde{T}/(\epsilon^{1/3})\rceil$. Then it is clear that $t_{m_J}\leq 24\tilde{T}$.
Thus Proposition \ref{pro4-1} and the domain Markov property of harmonic
explorer process imply that there exists a constant $c>0$ depending only on $K$ such that
\begin{equation}\label{eq12}
\mid\mathbb{E}[W^\epsilon(t_{m_{j+1}})-W^\epsilon(t_{m_{j}})\mid\mathcal{F}_{j}]\mid\leq{c\epsilon^{1/2}}
\end{equation}
and
\begin{equation}\label{eq13}\mid\mathbb{E}[(W^\epsilon(t_{m_{j+1}})-W^\epsilon(t_{m_{j}}))^{2}-4(t_{m_{j+1}}
-t_{m_{j}})\mid\mathcal{F}_{j}]\mid\leq{c\epsilon^{1/2}}
\end{equation}
for $j=0,1,\dots,J-1$, where $\mathcal{F}_{j}$ is the $\sigma$-algebra generated by $\gamma^\epsilon[0,m_{j}]$.

Let
$$X_{j}:=(W^\epsilon(t_{m_{j}})-W^\epsilon(t_{m_{j-1}}))-\mathbb{E}[W^\epsilon(t_{m_{j}})-W^\epsilon(t_{m_{j-1}})\mid\mathcal{F}_{j-1}]$$
for $j=1,2,\dots,J$. We define a process $M$ by $M_{0}=0$ and $ M_{j}=\sum_{k=1}^{j}X_{k}$ for $j=1,2,\dots,J$.
Then it is clear that $M$ is a martingale with respect to $\mathcal{F}_{j}$.
Moreover, it follows from the definitions of $X_j$ and $m_j$ that
$$\parallel{M_{j}-M_{j-1}}\parallel_{\infty}=\|X_{j}\|_{\infty}\leq{4\epsilon^{1/6}}.$$
Thus applying Lemma \ref{lem5-2} to $M$ gives that there exist stopping times ${\tau_{j}}$ for the
standard Brownian motion $B$ and a coupling of $B$ with $M$
such that $M_{j}=B(\tau_{j})$ for $j=0,1,\dots,J$.

Next, we will estimate $|4t_{m_j}-\tau_j|$. By the triangle inequality, it suffices to estimate
$|4t_{m_j}-Y_{j}|$ and $|Y_{j}-\tau_j|$ respectively, where $Y_{j}:=\sum_{k=1}^{j}X_{k}^{2}$
denotes the natural time associated to $M$.
For the first term $|4t_{m_j}-Y_{j}|$ we have
\begin{eqnarray}\label{eq14-1}
\mathbb{P}\Big(\max_{1\leq{j}\leq{J}}\mid{Y_{j}-4t_{m_{j}}}
\mid&\geq&3\epsilon^{1/6}\mid\log\epsilon\mid\Big)=\mathbb{P}\Big(\max_{1\leq{j}\leq J}\mid\sum_{k=1}^{j}(X_{k}^{2}-Z_{k})
\mid\geq{3\epsilon^{1/6}|\log \epsilon|}\Big)\nonumber\\
&\leq&\mathbb{P}\Big(\max_{1\leq{j}\leq J}\mid\sum_{k=1}^{j}(X_{k}^{2}-\mathbb{E}[X_{k}^{2}\mid\mathcal{F}_{k-1}])
\mid\geq\epsilon^{1/6}|\log \epsilon|\Big)\nonumber\\
&&+\mathbb{P}\Big(\max_{1\leq{j}\leq J}\mid\sum_{k=1}^{j}(\mathbb{E}[X_{k}^{2}
\mid\mathcal{F}_{k-1}]-\mathbb{E}[Z_{k}\mid\mathcal{F}_{k-1}])\mid\geq\epsilon^{1/6}|\log \epsilon|\Big)\nonumber\\
&&+\mathbb{P}\Big(\max_{1\leq{j}\leq J}\mid\sum_{k=1}^{k}(Z_{k}-\mathbb{E}[Z_{k}\mid\mathcal{F}_{k-1}])\mid
\geq\epsilon^{1/6}|\log \epsilon|\Big)\nonumber\\
&\triangleq&{A_{1}+A_{2}+A_{3}},
\end{eqnarray}
where $Z_{j}=4t_{m_{j}}-4t_{m_{j-1}}$. Applying Lemma \ref{lem5-3} to $A_{1}$
with $a=\epsilon^{1/6}\mid\log\epsilon\mid, b=\epsilon^{1/6}, c=e^{-2}\epsilon^{1/3}|\log\epsilon|$, combined with
definitions of $X_j$ and $J$,  implies that $A_1=O(\epsilon^{1/6})$ for $\epsilon$ small enough.
From (\ref{eq12}) and (\ref{eq13}) we get that $A_{2}=0$ if $\epsilon$ is small enough.
By the same argument of estimating $A_{1}$ and together with the inequality $\max\mid{Z_{j}}\mid\leq 8\epsilon^{1/3}$ we deduce
that $A_{3}=O(\epsilon^{1/6})$. Hence we conclude from (\ref{eq14-1}) that
\begin{equation}\label{eq14}\mathbb{P}\Big(\max_{1\leq{j}\leq{J}}\mid{Y_{j}-4t_{m_{j}}}
\mid\geq3\epsilon^{1/6}\mid\log\epsilon\mid\Big)=O{(\epsilon^{1/6})}
\end{equation}
for  all $\epsilon$ small enough.

For the second term $|Y_{j}-\tau_{j}|$, let $\tilde{Z}_{j}=\tau_{j}-\tau_{j-1}$ and
let $\mathcal{G}_{j}$  denote the $\sigma$-algebra generated by $B[0, \tau_{j}]$.
Then we get
\begin{eqnarray}\label{eq16-1}
\mathbb{P}\Big(\max_{1\leq{j}\leq{J}}
\mid{Y_{j}-\tau_{j}}\mid&\geq&3\epsilon^{1/6}\mid\log\epsilon\mid\Big)=\mathbb{P}\Big(\max_{1\leq{j}\leq J}\mid\sum_{k=1}^{j}(X_{k}^{2}-\tilde{Z}_{k})\mid\geq{3\epsilon^{1/6}|\log\epsilon|}\Big)\nonumber\\
&\leq&\mathbb{P}\Big(\max_{1\leq{j}\leq J}\mid\sum_{k=1}^{j}(X_{k}^{2}-\mathbb{E}[X_{k}^{2}\mid\mathcal{G}_{k-1}])
\mid\geq\epsilon^{1/6}\mid\log\epsilon\mid\Big)\nonumber\\
&&+\mathbb{P}\Big(\max_{1\leq{j}\leq J}\mid\sum_{k=1}^{j}(\mathbb{E}[X_{k}^{2}
\mid\mathcal{G}_{k-1}]-\mathbb{E}[\tilde{Z}_{k}\mid\mathcal{G}_{k-1}])\mid\geq\epsilon^{1/6}|\log\epsilon|\Big)\nonumber\\
&&+\mathbb{P}\Big(\max_{1\leq{j}\leq J}\mid\sum_{k=1}^{j}(\tilde{Z}_{k}-
\mathbb{E}[\tilde{Z}_{k}\mid\mathcal{G}_{k-1}])\mid\geq\epsilon^{1/6}|\log\epsilon|\Big)\nonumber\\
&\triangleq&{\tilde{A}_{1}+\tilde{A}_{2}+\tilde{A}_{3}}.
\end{eqnarray}
The estimation of $\tilde{A}_{1}$ is identical to that of $A_{1}$ above.
Note that $X_{k+1}^{2}=(B_{\tau_{k+1}}-B_{\tau_{k}})^{2}$,
so we obtain from (\ref{eq9-1}) that $\tilde{A}_{2}=0$.
Applying Lemma \ref{lem5-3}  to $\tilde{A}_{3}$, together with (\ref{eq9-1}) and (\ref{eq9-2}), implies
that $\tilde{A}_{3}=O{(\epsilon^{1/6})}$. Thus we deduce from (\ref{eq16-1}) that
\begin{equation}\label{eq16}\mathbb{P}\Big(\max_{1\leq{j}\leq{J}}
\mid{Y_{j}-\tau_{j}}\mid\geq3\epsilon^{1/6}\mid\log\epsilon\mid\Big)=O{(\epsilon^{1/6})}
\end{equation}
for $\epsilon$ small enough.
Therefore,  from (\ref{eq14}) and (\ref{eq16}) we get
\begin{equation}\label{eq17}\mathbb{P}\Big(\max_{1\leq{j}\leq{J}}\mid{4t_{m_{j}}
-\tau_{j}}\mid\geq 6\epsilon^{1/6}\mid\log\epsilon\mid\Big)=O{(\epsilon^{1/6})}
\end{equation}
for $\epsilon$ sufficiently small.

At the same time, it follows from (\ref{eq9-3}) that
\begin{equation}\label{eq18}\sup\{\mid{B(t)-B(\tau_{j-1}})\mid:
t\in[\tau_{j-1},\tau_{j}]\}\leq4\epsilon^{1/6}
\end{equation}
for every $j\leq{J}$. By the definition of $m_{j}$ and (\ref{eq6}) one obtains
\begin{equation}\label{eq18-1}\sup\{\mid{W^\epsilon(t_{m_{j}})-W^\epsilon(t)}\mid: t\in[t_{m_{j-1}},t_{m_{j}}]\}\leq c\epsilon^{1/6},
\end{equation}
where $c>0$ is a constant.
The definitions of $M_j$ and $J$, combined with (\ref{eq12}), give
\begin{equation}\label{eq18-2}\sup\{\mid{W^\epsilon(t_{m_{j}})-M_{j}}\mid: j\leq J\}\leq{c\tilde{T}\epsilon^{1/6}}.
\end{equation}
The definitions of $Y_j$ and $t_{m_j}$ imply $(Y_{j+1}-Y_{j})+(t_{m_{j+1}}-t_{m_j})\geq \epsilon^{1/3}$. Summing over $j$
gives $Y_J+t_{m_J}\geq J\epsilon^{1/3}\geq 12\tilde{T}$.
This implies that the event that $t_{m_J}<2\tilde{T}$ is contained in the event
that  $\mid Y_J-4t_{m_J}\mid\geq 2\tilde{T}$.
It follows from (\ref{eq14}) that
\begin{equation}\label{eq19}\mathbb{P}(t_{m_J}<2\tilde{T})=O(\epsilon^{1/6}).\end{equation}

Lastly, we consider the following event
\begin{eqnarray*}E=\{t_{m_{J}}\geq2\tilde{T}\}\cap&\Big\{&\sup_{t\in[0, 2T-\lambda]}\sup_{s\in(0, \lambda]}\mid{B(t+s)-B(t)}\mid\leq{\sqrt{6\lambda\mid\log{\lambda}\mid}}\Big\}\\
\cap&\Big\{&\max_{j\leq{J}}\mid{4t_{m_{j}}-\tau_{j}}\mid\leq6\lambda\Big\},
\end{eqnarray*}
where $\lambda=\lambda(\epsilon)=\epsilon^{1/6}\mid\log\epsilon\mid$.
Then from (\ref{eq17}),(\ref{eq19}) and Lemma \ref{lem5-4} (with $\delta=1$, and $\mu=\sqrt{6\mid\log{\lambda}\mid}$)
we get
\[\mathbb{P}(E^{c})=O(\epsilon^{1/6}\mid\log\epsilon\mid).\]

On the event $E$, by the triangle inequality we have
\begin{eqnarray*}
\sup\{\mid W^\epsilon(t)&-&B(4t)\mid: t\in[0, \tilde{T}]\}
\nonumber\\&\leq&\max_{1\leq{j}\leq{J}}\Big(\sup\{\mid(W^\epsilon(t)-W^\epsilon(t_{m_j})\mid: t\in[t_{m_{j-1}},t_{m_{j}}]\}
\nonumber\\&&+\mid W^\epsilon(t_{m_{j}})-B(\tau_{j})\mid+\sup\{\mid{B(\tau_{j})-B(4t)}\mid: t\in[t_{m_{j-1}},t_{m_{j}}]\}\Big).
\end{eqnarray*}
It follows from (\ref{eq18-1}) and (\ref{eq18-2}) that the first two terms on the right hand of the above
inequality are $O(\tilde{T}\epsilon^{1/6})$ uniformly in $j$. For the last term, we obtain from
(\ref{eq18}) that
\begin{eqnarray*}
\sup\{\mid B(\tau_{j})&-&B(4t)\mid: t\in[t_{m_{j-1}},t_{m_{j}}]\}
\nonumber\\&=&\sup\{\mid B(\tau_{j})-B(s)\mid: s\in[4t_{m_{j-1}},4t_{m_{j}}]\}
\nonumber\\&\leq&\sup\{\mid B(\tau_{j})-B(s)\mid: s\in[\tau_{j-1}-6\lambda, \tau_{j}+6\lambda]\}
\nonumber\\&\leq&4\epsilon^{1/6}+\sup\{\mid B(\tau_{k-1})-B(s)\mid: s\in[\tau_{k-1}-6\lambda,\tau_{k-1}]\}
\nonumber\\&&+\sup\{\mid B(\tau_{j})-B(s)\mid: s\in[\tau_{j},\tau_{j}+6\lambda]\}
\nonumber\\&\leq&4\epsilon^{1/6}+c(\lambda|\log\lambda|)^{1/2}.
\end{eqnarray*}
Thus we can couple $W^\epsilon(t)$ and $B(t)$ so that
\begin{equation}\label{eq19-1}\mathbb{P}\Big(\sup_{t\in[0,\tilde{T}]}\{\mid W^\epsilon(t)-B(4t)\mid\}>c_{1}
\tilde{T}\epsilon^{1/12}|\log\epsilon\log\lambda|^{1/2}\Big)<c_{2}\epsilon^{1/6}\mid\log\epsilon\mid,
\end{equation}
where $c_1>0,c_2>0$ are universal constants. A direct computation gives that for any fixed $0<\nu<1/12$,
$\epsilon^{1/12-\nu}>c_{1}\tilde{T}\epsilon^{1/12}|\log\epsilon\log\lambda|^{1/2}$ for $\epsilon$ sufficiently small.
So we deduce from (\ref{eq19-1}) that (\ref{eq1-1}) holds for $t\in[0,\tilde{T}]$ if $\epsilon$ is small enough.
Note that the scaled Brownian motion $B(4t)$ is unlikely to hit $\{-R,R\}$ before time
$24T$ if $R$ is sufficiently large. Moreover, the definition of $\tilde{T}$ gives that $\tilde{T}\rightarrow T$
when $R\rightarrow\infty$. Hence we conclude that (\ref{eq1-1}) also holds for $t\in [0,T]$ by
taking  limit as $R\rightarrow\infty$. Thus we finish the proof of the theorem.  $\hfill\square$

\section{Convergence rate for harmonic explorer path}\label{ecrhep}

In \cite{SS05} Schramm and Sheffield proved that as $\epsilon\rightarrow 0$, the image of the harmonic
explorer path $\gamma^\epsilon$ in $\mathbb{H}$, $\tilde{\gamma}^\epsilon$, converges uniformly with respect to a natural
metric on curves to the chordal $\mbox{SLE}_4$ path starting from $0$; see \cite[Theorem 3.3]{SS05}
for a precise statement. The goal of this section is to prove Theorem \ref{thm2}, which may be viewed as a quantitative
version of \cite[Theorem 3.3]{SS05}. To this end, we first establish the corresponding estimate of Loewner curves
in the deterministic setting; the details are given in Proposition \ref{pro6-2}. Next, we will show
that the assumptions of Proposition \ref{pro6-2} are satisfied on an event of large probability.

\subsection{The estimate for deterministic Loewner curves}

In this subsection we will consider a deterministic setting with two solutions to (\ref{eq02}) driven by
functions which are at uniform distance at most $\epsilon>0$.
If the growth of the derivative of one solution is known and the Loewner curve corresponding to the other
solution satisfies the John-type condition,  the supremum distance between the corresponding two curves can be controlled.
To be more precise, we provide the following proposition, which is analogous to \cite[Lemma 3.4]{joh}
concerning the radial case.

\begin{prop}\label{pro6-2} Let $(f_j,U_j)$ be $\mathbb{H}$-Loewner pairs generated by the curves $\gamma_j$ for $j=1,2$.
Fix $T<\infty$. Assume that there exist $\epsilon>0, \beta<1$ and $p,r\in(0,1)$ such that the
following holds.

(a) The driving terms satisfy
\[\sup\limits_{t\in[0,T]}|U_1(t)-U_2(t)|\leq \epsilon;\]

(b) There exists a constant $c<\infty$ such that the derivative estimate
\[\sup\limits_{t\in[0,T]}|f_1'(t,U_1(t)+id)|\leq cd^{-\beta}\]
holds for any $d\leq \epsilon^p$;

(c) There exists a constant $\tilde{c}<\infty$ such that the tip structure modulus for $(\gamma_2(t),t\in[0,T])$ in $\mathbb{H}$
satisfies $\eta(\epsilon^p)\leq \tilde{c}\epsilon^{pr}$;

Then there is a constant $\hat{c}=\hat{c}(T,\beta,p,r,c,\tilde{c})<\infty$ such that
\begin{equation}\label{equ6-2}
\sup\limits_{t\in[0,T]}|\gamma_1(t)-\gamma_2(t)|\leq \hat{c}\max\{\epsilon^{p(1-\beta)r},\epsilon^{(1-p)r}\}.
\end{equation}
\end{prop}

To prove Proposition \ref{pro6-2}, we need the following lemmas. First, using the reverse chordal Loewner equation
it is not hard to derive the following estimate of the uniform distance between two solutions to (\ref{eq02}) with the
driving terms whose supremum distance is at most $\epsilon$;
see \cite[Lemma 2.3]{jrw} for a proof.

\begin{lemma}\label{lem6-3} Let $(f_j,U_j)$ be $\mathbb{H}$-Loewner pairs for $j=1,2$. Fix $T<\infty$ and suppose that
\[\sup\limits_{t\in[0,T]}|U_1(t)-U_2(t)|\leq \epsilon\]
for arbitrary $\epsilon>0$. Then there exists a constant $c$ depending only on $T$ such that
for any $z_1,z_2\in\mathbb{H}$ with $|z_1-z_2|\leq \epsilon$ and $\mbox{Im}z_j\geq d,j=1,2$,
\[|f_1(t,z_1)-f_2(t,z_2)|\leq c\epsilon d^{-1}\]
for each $t\in(0,T]$.
\end{lemma}

Next, for a chordal Loewner curve we have the following result that the distance between the tip
of curve and the corresponding solution can be bounded by  the tip structure modulus up to
a multiplicative constant.

\begin{lemma}\label{lem6-4} For a given $T<\infty$, there exist constants $0<c_1,c_2,c_3<\infty$ with $c_1,c_2$ depending only
on $T$ and $c_3$ universal such that the following holds.  Let $\gamma$ be a curve in $\mathbb{H}$ generating the Loewner
pair $(f,U)$ and let $\eta(\delta)$ be the tip structure modulus for $(\gamma(t),t\in[0,T])$. Then if $t\in[0,T]$
and $\Lambda_{f(t)}(d)<c_1$ where $\Lambda_{f(t)}(d):=\mbox{dist}(f(t,U_t+id),\partial H_t)$,  one has
\begin{equation}\label{equ6-2-02}|\gamma(t)-f(t,U_t+id)|\leq c_2\eta(c_3\Lambda_{f(t)}(d))
\end{equation}
for any small $d>0$.
\end{lemma}

\noindent\textbf{Proof.} The lemma is similar to \cite[Proposition 3.2]{joh} on radial Loewner curves.
Observe that the main difference between
the tip structure modulus for a radial Loewner curve $\tilde{\gamma}$ and the one for a chordal Loewner curve $\gamma$ is that
crosscuts separate $\tilde{\gamma}(t)$ from $0$ in the former, while crosscuts separate $\gamma(t)$ from $\infty$
in the latter. The remaining constructions are the same. So similarly to the proofs of
\cite[Lemma 3.1 and Proposition 3.2]{joh}
we can deduce that there exist  universal constants $\tilde{c}_2>0$ and $\tilde{c}_3>0$ such that
\[\rho(\gamma(t),f(t,U_t+id))\leq \tilde{c}_2\eta(\tilde{c}_3\Lambda_{f(t)}(d)),\]
where the metric $\rho(\cdot,\cdot)$ is defined in Section \ref{haex}. This, combined with (\ref{ineq1}),
implies that (\ref{equ6-2-02}) holds. $\hfill\square$

\

\noindent\textbf{Proof of Proposition \ref{pro6-2}.} By the triangle inequality we have
\begin{eqnarray*}
|\gamma_1(t)-\gamma_2(t)|&\leq&|\gamma_1(t)-f_1(t,U_1(t)+i\epsilon^p)|
+|f_1(t,U_1(t)+i\epsilon^p)-f_1(t,U_2(t)+i\epsilon^p)|\\
&&+|f_1(t,U_2(t)+i\epsilon^p)-f_2(t,U_2(t)+i\epsilon^p)|
+|f_2(t,U_2(t)+i\epsilon^p)-\gamma_2(t)|\\
&\triangleq & b_1+b_2+b_3+b_4.
\end{eqnarray*}
From the condition (b) it is easy to see  that $b_1\leq c\epsilon^{p(1-\beta)}$. Since  $\epsilon^p\geq\epsilon$, by
 Koebe's distortion theorem and  the condition (b) we can deduce that $b_2\leq c\epsilon^{p(1-\beta)}$,
where the constant $c>0$ may change from line to line.
In addition, from Lemma \ref{lem6-3} we obtain that
$b_3\leq c\epsilon^{(1-p)}$,
where $c>0$ is a constant depending only on $T$.

It remains to estimate $b_4$. We claim that $\Lambda_{f_2(t)}(\epsilon^p)\leq c(\epsilon^{p(1-\beta)}+\epsilon^{1-p})$,
where $\Lambda_{f_2(t)}(\cdot)$ is defined  in Lemma \ref{lem6-4}.
Indeed, by Cauchy's inequality and Lemma \ref{lem6-3} we get
\begin{equation}\label{equ6-2-5}|\epsilon^p|f_1'(t,U_2(t)+i\epsilon^p)|-\epsilon^p|f_2'(t,U_2(t)+i\epsilon^p)||\leq
2\max_{z\in D_{\epsilon^p/2}}|f_1(t,z)-f_2(t,z)|\leq 4c\epsilon^{(1-p)},
\end{equation}
where $D_{\epsilon^p/2}:=\{z\in\mathbb{H}:|z-(U_2(t)+i\epsilon^p)|\leq \epsilon^p/2\}$.
By  Koebe's estimate and the triangle inequality we have
\begin{eqnarray*}\Lambda_{f_2(t)}(\epsilon^p)&\leq& 4\epsilon^p|f_2'(t,U_2(t)+i\epsilon^p)|\\
&\leq& 4|\epsilon^p|f_1'(t,U_2(t)+i\epsilon^p)|+4|\epsilon^p|f_1'(t,U_2(t)+i\epsilon^p)|-\epsilon^p|f_2'(t,U_2(t)+i\epsilon^p)||,
\end{eqnarray*}
which, combined with the condition (b) and (\ref{equ6-2-5}), implies that the claim holds.
Thus from Lemma \ref{lem6-4}, the condition (c) and the claim  it follows that
$b_4\leq c(\epsilon^{p(1-\beta)r}+\epsilon^{(1-p)r})$.
Hence we conclude that (\ref{equ6-2}) holds. $\hfill\square$

\subsection{The estimate for harmonic explorer paths}
In this subsection we will show that  the conditions (b) and (c) in Proposition \ref{pro6-2} are
satisfied with a large probability, which concludes Theorem \ref{thm2}.

To accomplish this, the following results are needed. First, the derivative estimate for chordal $\mbox{SLE}_4$
will be applied. Let
\[\chi(\beta)=\frac{3}{2}\beta+\frac{1}{2}+\frac{\beta^2}{2(1+\beta)}\]
and
\[q(\beta)=\min\{\frac{17}{8}\beta,\chi(\beta)-2\}\]
for $\sqrt{3}/2<\beta<1$. Then we have the following estimate on the growth of the derivative for chordal $\mbox{SLE}_4$;
see \cite[Proposition 4.2]{jl} and \cite[Proposition A.4]{joh} for a proof.

\begin{lemma}\label{lem6-5} Let $T<\infty$ be fixed and let $(f(t,z),B(4t))$ be the chordal $\mbox{SLE}_4$ Loewner pair.
Let $\beta\in(\sqrt{3}/2,1)$ and $q<q(\beta)$. Then exists a constant $0<c<\infty$ depending
only on $T$ and $q$ such that for each $y_*<1$,
\[\mathbb{P}\{\forall y\leq y_*,\sup\limits_{t\in [0,T]} |f'(t,iy+B(4t))|\leq cy^{-\beta}\}\geq 1-cy_*^q.\]
\end{lemma}

Next, we will need an estimate for the tip structure modulus of harmonic explorer path, which is given as follows.

\begin{prop}\label{pro6-6} Suppose that $D^\epsilon$ is a $TG^\epsilon$ domain approximation of $D$ with two
prescribed points $\hat{v}_0,\hat{v}_e\in \partial D^\epsilon$. Let
$\gamma^\epsilon$ be the harmonic explorer path from $\hat{v}_0$ to $\hat{v}_e$ in $D^\epsilon$ defined in
Section \ref{haex}. Let $\eta(\delta)$ be the tip structure
modulus for $\gamma^\epsilon$ stopped when first reaching distance $\sigma>0$ from $\hat{v}_e$.  Set $r\in(0,1)$.
Then there
exist a universal constant $c_0>0$ and $c=c(r,\sigma)<\infty$ such that
\begin{equation}\label{eq6-1}\mathbb{P}\{\eta(\delta)\leq\delta^r\}\geq 1-c\delta^{1/2-r/2}
\end{equation}
for $\epsilon$ sufficiently small and $\delta>c_0\epsilon$.
\end{prop}

To prove Proposition \ref{pro6-6}, we need two more lemmas. The first one is the following Beurling's estimate
for a random walk; see \cite{L05} or \cite{ll} for a proof.

\begin{lemma}\label{lem6-7} (Beurling's estimate) There exists a constant $c<\infty$ such that the following holds.
 Let $\Omega\subset TG^\epsilon$ be an
infinite connected set. Let $S$ be a simple random walk on $TG^\epsilon$ started from $z$ and stopped at the first time $\tau_\Omega$
at which $S$ hits $\Omega$. Then for $r>1$, one has
\[\mathbb{P}\{|S(\tau_\Omega)-z|\geq r\mbox{dist}(z,\Omega)\}\leq cr^{-1/2}.\]
\end{lemma}
The next one is the revisiting probability estimate for the harmonic explorer path $\gamma^\epsilon$.

\begin{lemma}\label{lem6-8} Let $0<r<R<\infty$. Let $B(z,r)$ be a ball of radius $r$ intersecting $\gamma^\epsilon[0,n]$,
and let $B(z,R)$ be the concentric ball with radius $R$.
Suppose that there is no component of $B(z,R)\cap(D^\epsilon\setminus
\gamma^\epsilon[0,n])$ whose boundary intersects both $\partial_+\hat{D}_n^\epsilon$ and $\partial_-\hat{D}_n^\epsilon$,
where $\partial_+\hat{D}_n^\epsilon$ and $\partial_-\hat{D}_n^\epsilon$ are defined  in Section \ref{rcmo}. Then one has
\begin{equation}\label{eq6-2}\mathbb{P}[\gamma^\epsilon[n,N]\cap B(z,r)\neq\emptyset|\gamma^\epsilon[0,n]]\leq O(1)(r/R)^{1/2}.
\end{equation}
\end{lemma}

\noindent\textbf{Proof.} This is essentially \cite[Proposition 6.3]{SS05}; the only difference is that the exponent $\hat{c}$ in the
right-hand side of (\ref{eq6-2}) was not specified in \cite{SS05}. So it is enough to explain how one gets the exponent
$1/2$. Indeed, in the proof of \cite[Proposition 6.3]{SS05}, that is, on p.17, before the inequality (6.4), we apply Lemma \ref{lem6-7} to
see that $\hat{c}=1/2$. This, combined with the Markovian property for the harmonic explorer path $\gamma^\epsilon$, implies
the lemma. $\hfill\square$

\

\noindent\textbf{Proof of Proposition \ref{pro6-6}.} For $0<\delta\leq R$, the definition of tip structure modulus
implies that the event $\{\eta(\delta)\leq R\}$ is equivalent to the event $E$ that
$\gamma^\epsilon(t),t\in[0,T]$ has no nested $(\delta,R)$-bottleneck in $D^\epsilon$. In the following we will estimate the
probability of $E$.

We claim that $E^c$ is contained in the event $\tilde{E}$ that
there is a ball $B(z,\delta)$ of radius $\delta$ such that $\gamma[n,N]\cap B(z,\delta)\neq\emptyset$ conditioned on the
ball $B(z,\delta)$ intersecting $\gamma^\epsilon[0,n]$, where $0\leq n<N$.
Indeed, it is easy to see that
$E^c$ is contained in the event $\tilde{E}_1$ that there exists one
nested $(\delta,R)$-bottleneck in $D^\epsilon$ for $\gamma^\epsilon$ stopped when reaching $\partial B(\hat{v}_e,\sigma)$.
If $\tilde{E}_1$ occurs,  there must be some $t\in[0,T]$ and a crosscut $C_{t,\delta}$ of $D_t^\epsilon$ with diameter at
most $\delta$ that separates $\gamma^\epsilon(t)$ from $\hat{v}_e$ in $D_t^\epsilon$ such that
$\text{diam}(\gamma^\epsilon[t-\tau,t])\geq R$, where $\tau=\inf\{s>0:\gamma^\epsilon[t-s,t]\cap C_{t,\delta}\neq\emptyset\}$.
In this situation, there are three possibilities for the crosscut $C_{t,\delta}$: (i) $C_{t,\delta}$ intersects only the path $\gamma^\epsilon[0,t]$; (ii)
$C_{t,\delta}$ intersects both $\gamma^\epsilon[0,t]$ and the boundary $\partial D^\epsilon$;  and (iii) $C_{t,\delta}$
intersects only  $\partial D^\epsilon$. Notice that $D$ is fixed and that $D^\epsilon$ tends to $D$ as $\epsilon\rightarrow0$,
which implies that the case (iii) never occurs for $\epsilon$ and $\delta$ sufficiently small.
Since  the crosscut $C_{t,\delta}$ can be taken to be a subarc of $\partial B(z,\delta)$, it follows that $\tilde{E}$ must occur if $\tilde{E}_1$
happens with the case (i) or (ii). This implies that $\tilde{E}_1\subset\tilde{E}$. Hence we deduce that the claim holds.

Furthermore, from the definition of nested $(\delta,R)$-bottleneck for $\gamma^\epsilon$
we deduce that $B(z,R)\cap D_n^\epsilon$
has no component whose boundary intersects both $\partial_+\hat{D}_n^\epsilon$ and $\partial_-\hat{D}_n^\epsilon$,
where $B(z,R)$ denotes the concentric ball with radius $R$.
Thus applying Lemma \ref{lem6-8} to the event $\tilde{E}$ gives
\[\mathbb{P}[\tilde{E}]
\leq O(1)(\delta/R)^{1/2}.\]
This, combined with the claim above,  implies
\[\mathbb{P}[E]\geq 1-\mathbb{P}[\tilde{E}]\geq 1-O(1)\delta^{1/2-r/2},\]
if we take $R=\delta^r,r\in(0,1)$. Hence from the above equivalence we can conclude that
the inequality (\ref{eq6-1}) holds. $\hfill\square$

Lastly, we need the following result on the relationship between the tip structure modulus for a Loewner curve
and the one for the image of the curve under a conformal map, which will be proved in the Appendix.

\begin{lemma}\label{pro6-9} Suppose $D$ is a simply connected Jordan domain with $C^{1+\alpha}$ boundary, where $\alpha>0$. Let
$D^\epsilon$ be the $TG^\epsilon$ domain approximation of $D$ with two prescribed points $\hat{v}_0,\hat{v}_e\in \partial D^\epsilon$
and let $\gamma^\epsilon$ be a Loewner curve in $D^\epsilon$ from $\hat{v}_0$ to $\hat{v}_e$. There is a constant $c$ depending only
on $\alpha$ and the diameter of $D$ such that the following holds. Set $0<r<1/2$ and $d_\epsilon=\epsilon^r$.
Let $\eta(\delta;D^\epsilon)$ be the tip structure modulus for $\gamma^\epsilon$ in $D^\epsilon$. Then for
$\epsilon$ sufficiently small the tip structure modulus $\eta(\delta;\mathbb{H})$ for $\phi(\gamma^\epsilon)$
in $\mathbb{H}$ satisfies
\[\eta(c^{-1}d_\epsilon;\mathbb{H})\leq c\eta(d_\epsilon;D^\epsilon).\]
\end{lemma}

\noindent\textbf{Proof of Theorem \ref{thm2}.} It follows from Theorem \ref{thm1} that there exists a coupling of the chordal $\mbox{SLE}_4$
path $\tilde{\gamma}$ and the image of the harmonic explorer path $\tilde{\gamma}^\epsilon$. We will estimate the distance between $\tilde{\gamma}$
and $\tilde{\gamma}^\epsilon$ in this coupling. Take $\sigma\in(0,1/12)$ and $\epsilon<\epsilon_0$, where $\epsilon_0$ is defined in Theorem \ref{thm1};
for $p\in(0,1)$, set
\[\varrho_\epsilon=\epsilon^\sigma, d_\epsilon=(\varrho_\epsilon)^p.\]
For each $\epsilon\leq \epsilon_0$, we shall define three events each of which occurs with large probability
in our coupling. On the intersection of these events, we are able to apply those estimates in Proposition \ref{pro6-2}.

(i) Let $E_1^\epsilon$ denote the event that
\[\sup\limits_{t\in[0,T]}|W^\epsilon(t)-B(4t)|\leq \varrho_\epsilon.\]
 Then Theorem \ref{thm1} gives that there exists $\epsilon_0>0$ such that for any $\epsilon<\epsilon_0$,
\[\mathbb{P}(E_1^\epsilon)\geq 1-\varrho_\epsilon.\]

(ii) For $\beta\in(\sqrt{3}/2,1)$, let $E_2^\epsilon=E_2(\sigma,\beta,T,c)$
be the event that the chordal $\mbox{SLE}_4$ Loewner chain $(f(t,\cdot))$ driven  by $B(4t)$ satisfies the estimate
\[\sup\limits_{t\in[0,T]}|f'(t,B(4t)+id)|\leq c d^{-\beta}\]
for each $d\leq d_\epsilon$. Then by Lemma \ref{lem6-5} there exist $c'<\infty$, independent of $\epsilon$ and $\epsilon_1>0$,
such that
\[\mathbb{P}(E_2^\epsilon)\geq 1-c'(d_\epsilon)^q\]
for every $\epsilon<\epsilon_1$, where
\[q<q(\beta)=-\frac{3}{2}+\frac{3\beta}{2}+\frac{\beta^2}{2(1+\beta)}.\]

(iii) For $r\in(0,1)$, let $E_3^\epsilon=E_3(\sigma,r,p,\tilde{c},\mbox{diam}(D))$ be the event that the tip
structure modulus $\eta$ for $\tilde{\gamma}^\epsilon(t)$ with $t\in[0,T]$ in $\mathbb{H}$ satisfies
\[\eta(d_\epsilon)\leq \tilde{c}(d_\epsilon)^r.\]
Then by Lemma \ref{pro6-9} and Proposition \ref{pro6-6}, we know that there exists $\tilde{c}'<\infty$,
independent of $\epsilon$ and $\epsilon_2>0$, such that
if $\epsilon<\epsilon_2$ then
\[\mathbb{P}(E_3^\epsilon)\geq 1-\tilde{c}'(d_\epsilon)^{1/2-r/2}.\]

Consequently, there  is a constant $c>0$ depending only on
$\sigma,r,p,T,\beta$ and $\mbox{diam}(D)$ such that for all $\epsilon$ sufficiently small,
\begin{equation}\label{eq6-3-1}\mathbb{P}(E_1^\epsilon\cap E_2^\epsilon\cap E_3^\epsilon)\geq 1-c(\varrho_\epsilon+(d_\epsilon)^q
+(d_\epsilon)^{1/2-r/2}).
\end{equation}
The benefit is that on the event $E_1^\epsilon\cap E_2^\epsilon\cap E_3^\epsilon$, applying Proposition \ref{pro6-2} with $c=c',\tilde{c}=\tilde{c}'$ independent of
$\epsilon$, we can obtain that there exists a constant $\hat{c}>0$ independent of $\epsilon$ (but dependent on the above parameters)
such that for all $\epsilon$ sufficiently small,
\begin{equation}\label{eq6-3}\sup\limits_{t\in[0,T]}|\tilde{\gamma}^\epsilon(t)-\tilde{\gamma}(t)|
\leq \hat{c}((d_\epsilon)^{r(1-\beta)}+(\varrho_\epsilon)^{r(1-p)}).
\end{equation}

We now intend to optimize  the exponents. Since $d_\epsilon=(\varrho_\epsilon)^p$, it is easy
to see that  $(d_\epsilon)^{r(1-\beta)}$ dominates in (\ref{eq6-3}) when $p\in(0,1/(2-\beta))$,
and $(\varrho_\epsilon)^{r(1-p)}$ dominates whenever $p\in[1/(2-\beta),1]$.
Meanwhile, observe that as $p$ tends to $1/(2-\beta)$, the convergence exponents obtained are the same for these two cases.
So it suffices to consider the first case $p\in(0,1/(2-\beta))$.
We define
\[\nu(\beta,r):=\min\Big\{r(1-\beta),-\frac{3}{2}+\frac{3\beta}{2}+\frac{\beta^2}{2(1+\beta)},\frac{1}{2}-\frac{r}{2}\big\}.\]
Then the optimal rate is obtained by first optimizing over $\beta,r$ and next choosing $p$ very close to $1/(2-\beta)$.
Let $\beta_*=0.9079\in
(\sqrt{3}/2,1)$ be a solution to
\[8\beta^3-14\beta^2-6\beta+11=0\]
and take $r_*=1/(3-2\beta_*)=0.8445\in (0,1)$. Then we have
\[\nu(\beta_*,r_*)=\max\{\nu(\beta,r):\sqrt{3}/2<\beta<1,0<r<1\}=0.0778.\]
Consequently, for any
\[m<m_*=\frac{\nu(\beta_*,r_*)}{2-\beta_*}=0.0712,\]
we obtain bounds in (\ref{eq6-3-1}) and (\ref{eq6-3}) of order $(\varrho_\epsilon)^m$ for all $\epsilon$ sufficiently small.
Since $1/15<m_*$, this implies that (\ref{eq1-2}) holds for $\tilde{\gamma}^\epsilon(t)$ and $\tilde{\gamma}(t)$,
where $t\in[0,T]$.

Moreover, from the proof of Lemma \ref{pro6-9} we get that there is a constant $c>0$ depending only on
$\alpha$ and the diameter of $D$ such that
\[c^{-1}|z-z'|\leq \rho(\phi(z),\phi(z'))\leq c|z-z'|\]
for $z,z'\in D^\epsilon$. This, combined with (\ref{ineq1}), implies
\[\tilde{c}^{-1}|\gamma^\epsilon(t)-\gamma(t)|\leq |\tilde{\gamma}^\epsilon(t)-\tilde{\gamma}(t)|
\leq \tilde{c}|\gamma^\epsilon(t)-\gamma(t)|\]
for $t\in[0,T]$, where $\tilde{c}>0$ is a constant depending only on $c,T$ and the diameter of $D$.
Hence we conclude that (\ref{eq1-2}) also holds for $\gamma^\epsilon$ and $\gamma$,
which completes the proof of the theorem. $\hfill\square$

\section{Appendix:  bounds on tip structure modulus}\label{appen}
In this Appendix we will begin with some concepts involved in the proof of  Lemma \ref{pro6-9},
and then prove Lemma \ref{pro6-9}, which gives that the tip structure modulus of the image of a curve
under a conformal map can be bounded by that of the curve times a constant under the hypothesis that
the boundary of domain is sufficiently regular.

\subsection{Some concepts}
In this subsection we briefly review the definitions of boundary structure modulus and quasicircle, and give
a result concerning the relationship between boundary structure modulus and conformal maps; see \cite{pom,wa}
for more details.
\begin{definit}
Let $D\subset\mathbb{C}$ be a simply connected bounded domain which contains $0$.
Let $\gamma$ be a cross-cut of $D$ which does not pass through $0$, and let $D_1$ be the one of
the two subregions of $D$ which does not contain $0$. For any $\delta > 0$, by considering all possible cross-cuts
$\gamma$ of $D$ with $\mbox{diam}(\gamma)\leq\delta$, the structure modulus $\eta_b(\delta)$ of the boundary of $D$,
is defined by
\[\eta_b(\delta):=\sup\limits_{diam(\gamma)\leq \delta}\mbox{diam}(D_1),\]
where $\mbox{diam}(\cdot)$ denotes the diameter of a curve or a domain.
\end{definit}

The structure modulus $\eta_b(\delta)$ is in a certain sense a measure
of the irregularity of $\partial D$. For the relationship between structure modulus and conformal maps,
we have the following result, see \cite[Theorem VII]{wa} for a proof.

\begin{lemma}\label{lem2-5}\quad Suppose that $D, D'$ are two simply connected bounded domains
such that $D'\subset D$ and $z=0$ lies in $D'$. Let $\eta_b(\delta)$ denote the structure modulus of $\partial D'$
and $\varrho$  the distance from $0$ to $\partial D'$. Suppose that $\mbox{dist}(\partial D, \partial D')<\epsilon$,
$0 <\epsilon < 1, \epsilon <\varrho/64$. Let
$\phi$ and $\psi$ be the two conformal maps from $D$ and $D'$, respectively, onto $\mathbb{U}$
normalized by $\phi(0)=0,\phi'(0)>0$ and $\psi(0)=0,\psi'(0)>0$.
If $\eta_b(\delta)\leq c\delta$ for some constant $c>0$, then for any $z\in D'$,
\[|\phi(z)-\psi(z)|\leq C \epsilon^{1/2}\log\frac{2}{\epsilon},\]
where $C>0$ is a constant depending only on $c,\varrho$ and the diameter of $D$.
\end{lemma}

\begin{definit}
A Jordan curve $\Gamma$ is called a quasicircle if it is of bounded
turning; i.e., if there exists a constants $c>1$ such that for all
points $z_1,z_2\in \Gamma$,
\[\mbox{diam} (J(z_1,z_2))\leq c|z_1-z_2|,\] where
$J(z_1,z_2)$ is the subarc of $\Gamma$ connecting $z_1, z_2$ which has the smaller diameter.
\end{definit}

This geometric condition implies that a quasicircle cannot visit $z_1$, wander far away,
and then return to a point around $z_1$. The smallest $c$ above may be considered as a measure of
regularity. It is well known that every quasicircle is the image of a circle under a quasiconformal
map of the Riemann sphere.

\subsection{Proof of Lemma \ref{pro6-9}}
In this subsection we will prove  Lemma \ref{pro6-9}. To accomplish this,  we need the following two lemmas.

\begin{lemma}\label{lem6-10} Let $D\subset\mathbb{C}$ be a simply connected bounded domain
whose boundary is a quasicircle. Then the boundary $\partial D^\epsilon$
of $D^\epsilon$ is also a quasicircle which is  within distance $\epsilon$ from $\partial D$.
\end{lemma}

\noindent\textbf{Proof.} From the construction of $D^\epsilon$, it is easy to see that
each point of $\partial D^\epsilon$ is at distance less than $\epsilon$ from $\partial D$.
Thus it remains to show that $\partial D^\epsilon$ is a quasicircle. Let $z,w\in \partial D^\epsilon$
be any two points.

We first consider the case where $|z-w|<\epsilon$.
Since $\partial D^\epsilon$ is a Jordan curve which
is a subset of the edge set of triangulation $D^\epsilon$ with mesh $\epsilon$, we have that
\begin{equation}\label{eq6-2-1}
\mbox{diam}(J(z,w))\leq (2\sqrt{3}/3)|z-w|,
\end{equation}
where $J(z,w)\subset\partial D^\epsilon$ is the arc
with smaller diameter connecting $z$ and $w$.

Now assume that $|z-w|\geq \epsilon$. Let $\tilde{z}$ and $\tilde{w}$ be points
on $\partial D$ closest to $z$ and $w$, respectively. It is clear that $|z-\tilde{z}|$ and $|w-\tilde{w}|$ are both at most
$\epsilon$. Let $\gamma_1$ and $\gamma_2$ be the two line segments connecting $z$ with $\tilde{z}$
and $w$ with $\tilde{w}$ respectively. First, assume that the curve $\Gamma=J(z,w)\cup\gamma_1\cup\gamma_2$ separates
$J(\tilde{z},\tilde{w})$ from $z_0$ in $D$. Let $\triangle_j, j=1,2,\dots,N$, be those lattice triangles whose
faces are outside of $D^\epsilon$ but whose boundaries touch $J(z,w)$. By the construction of $D^\epsilon$ and
the Jordan curve theorem, since $\Gamma$ separates $z_0$ and $J(\tilde{z},\tilde{w})$, each $\triangle_j$ is intersected
by $\gamma_1\cup\gamma_2\cup J(\tilde{z},\tilde{w})$. This implies
\[\mbox{diam}(\Gamma) \leq \mbox{diam}(J(\tilde{z},\tilde{w}))+2\epsilon.\]
Since $\partial D$ is a quasicircle, there is a constant $A>0$ such that
\[\mbox{diam}(J(\tilde{z},\tilde{w}))\leq A |\tilde{z}-\tilde{w}|.\]
Hence, we get
\begin{equation}\label{eq6-2-2}\mbox{diam}( J(z,w))\leq \mbox{diam}(\Gamma)\leq A|z-w|+2(A+1)\epsilon.
\end{equation}

Next, if $\Gamma$ does not separate $J(\tilde{z},\tilde{w})$ from $z_0$ in $D$, then  $(\partial D^\epsilon\setminus J(z,w))\cup \gamma_1\cup\gamma_2$ does separate $J(\tilde{z},\tilde{w})$ from $z_0$ in $D$ since $\Gamma$ is a crosscut
of $D$. Thus, follwoing the same argument as in the previous paragraph, we obtain that
\[\mbox{diam}(\partial D^\epsilon\setminus J(z,w))\leq \mbox{diam} J(\tilde{z},\tilde{w})+2\epsilon.\]
Besides it follows from the definition of $J(z,w)$ that $\mbox{diam} J(z,w)\leq\mbox{diam}(\partial D^\epsilon\setminus
J(x,y))$. So we deduce that (\ref{eq6-2-2}) holds in this case too.

From (\ref{eq6-2-1}) and (\ref{eq6-2-2}), we conclude that
\[
\mbox{diam} J(z,w)\leq (3A+2)|z-w|
\]
for any two points $z,w\in\partial D^\epsilon$, which implies that $\partial D^\epsilon$ is a quasicircle.
So we finish the proof of the lemma. $\hfill\square$

\begin{lemma}\label{lem6-11} Let $D\subset\mathbb{C}$ be a simply connected bounded domain whose boundary is a quasicircle,
and fix $u_0,u_e\in \partial D$.  Let $D^\epsilon$ be the $TG^\epsilon$ domain approximation of $D$
which satisfies that $\hat{v}_0\rightarrow u_0,\hat{v}_e\rightarrow u_e$ as $\epsilon\rightarrow 0$. Let $\psi,\phi$ be two
conformal maps from $D$ and $D^\epsilon$, respectively, onto $\mathbb{H}$ with $\psi(u_0)=0,\psi(u_e)=\infty$
and $\phi(\hat{v}_0)=0,\phi(\hat{v}_e)=\infty$. Then there exists a constant $c>0$ depending
only on the measure of regularity of $\partial D$ and the diameter of $D$ such that
\begin{equation}\label{eq6-2-3}\sup_{z\in D^\epsilon}\rho(\psi(z),\phi(z))\leq c\epsilon^{1/2}|\log\epsilon|,
\end{equation}
where $\rho(\cdot,\cdot)$ is defined in Section \ref{haex}.
\end{lemma}

\noindent\textbf{Proof.} It follows from Lemma \ref{lem6-10}
that $\partial D^\epsilon$ is a quasicircle. This implies that there
exists a constant $c>0$ depending only on the measure of regularity of $\partial D$ such that
\begin{equation}\label{eq6-2-4}
\mbox{diam}(J(z,w))\leq c|z-w|,
\end{equation}
for any two points $z,w\in\partial D^\epsilon$, where $J(z,w)$ denotes
the subarc of $\partial D^\epsilon$
connecting $z$ and $w$ which has the smaller diameter.

By (\ref{eq6-2-4}) and the definition of boundary structure modulus, we deduce that there exists a constant $C$
depending only on $c$ and the diameter
of $D$ such that the structure modulus $\eta_b$ of $\partial D^\epsilon$ satisfies
\[\eta_b(\delta)\leq C\delta\]
for any small $\delta>0$.
Note that $D^\epsilon\subset D$ and Lemma \ref{lem6-10} gives that each point
on $\partial D^\epsilon$ is within distance $\epsilon$ of a point on $\partial D$.
So applying Lemma \ref{lem2-5} to $\varphi\circ\psi$ and $\varphi\circ\phi$
with a suitable  normalization, where $\varphi(z)=(z-i)/(z+i)$,
we obtain (\ref{eq6-2-3}) from the definition
of $\rho(\cdot,\cdot)$. $\hfill\square$

\

\noindent\textbf{Proof of Lemma \ref{pro6-9}.} Let $R_\epsilon=\eta(d_\epsilon;D^\epsilon)$. Without loss of generality
we can assume that $R_\epsilon\geq 2d_\epsilon$. From the definition of tip structure modulus it suffices
to show that there is a constant $c$ independent of $\epsilon$ such that for any annulus $A(z)
=\{w:d_\epsilon\leq |w-z|\leq R_\epsilon\}, z\in D^\epsilon$, we have
\begin{equation}\label{eq6-2-5}\phi(A(z)\cap D^\epsilon)\subset\{w:c^{-1}d_\epsilon\leq \rho(w,\phi(z))\leq c R_\epsilon\}\cap\mathbb{H}.
\end{equation}
Indeed, since $\partial D$ is $C^{1+\alpha}$, applying Kellogg's theorem \cite{gm} to the conformal map $\varphi\circ\psi:D\rightarrow\mathbb{U}$
implies that $\varphi\circ\psi$ and $(\varphi\circ\psi)^{-1}$ are in $C^{1+\alpha}(\overline{D})$ and $C^{1+\alpha}(\overline{\mathbb{U}})$,
respectively. In particular, $\varphi\circ\psi$ is bilipschitz on $\overline{D}$. This, combined with the definition
of $\rho(\cdot,\cdot)$, implies that there exists a constant $c>0$ depending only on
$\alpha$ and the diameter of $D$ such that
\begin{equation}\label{eq6-2-6}
c^{-1}|z-w|\leq \rho(\psi(z),\psi(w))|\leq c|z-w|
\end{equation}
for any $z,w\in \overline{D}$.

Meanwhile, the assumption on $\partial D$ implies that
$\partial D$ is a quasicircle. So for any $w=\phi(z')\in \phi(A(z)\cap D^\epsilon)$,
by the triangle inequality,  Lemma \ref{lem6-11} and
(\ref{eq6-2-6}), we can deduce that
\[c^{-1}d_\epsilon\leq \rho(\phi(z'),\phi(z))\leq c R_\epsilon,\]
which implies that (\ref{eq6-2-5}) holds. This completes the proof of the lemma. $\hfill\square$

\end{document}